\documentclass[twocolumn]{autart}    % Enable this line and disable the 
\usepackage{graphicx}
\usepackage[utf8]{inputenc}
\usepackage{mathptmx} % assumes new font selection scheme installed
\usepackage{newtxtext}
\usepackage{textcomp}
\usepackage[varg,bigdelims]{newtxmath}
\usepackage{enumitem}
\usepackage{dsfont, url}
\usepackage{mathtools}
\usepackage[usenames, dvipsnames]{xcolor}
\usepackage[%
colorlinks={true},%
citecolor={Blue},%
urlcolor={PineGreen},%
linkcolor={Sepia}%
]{hyperref}
\usepackage{float}
\usepackage{newfloat}
\usepackage[skip=1pt]{caption}

\usepackage{amsmath} % assumes amsmath package installed
\usepackage{amssymb}  % assumes amsmath package installed
\usepackage{todonotes}
\usepackage{pgf,tikz}
\usepackage{adjustbox}
\usepackage{verbatim} 
\usepackage{xargs}
\usepackage{cite}
\usetikzlibrary{plotmarks}
\usetikzlibrary{chains}
\usepackage{pgfplots}
\pgfplotsset{compat=newest}
\pgfplotsset{plot coordinates/math parser=false}
\newlength\figureheight
\newlength\figurewidth
\usetikzlibrary{shapes,arrows}
\usetikzlibrary{positioning, calc, fit}
\usetikzlibrary{decorations,decorations.pathmorphing,decorations.pathreplacing}
\usetikzlibrary{shapes.misc}
\usetikzlibrary{shadows,trees}
\pgfdeclaredecoration{switch}{initial}
{
	\state{initial}[width=+\pgfdecoratedinputsegmentremainingdistance]
	{
		\pgfsetlinewidth{1pt}
		\pgfpathcircle{\pgfpoint{0}{0}}{0.05\pgfdecoratedinputsegmentremainingdistance}
		\pgfpathlineto{\pgfpoint{0.25\pgfdecoratedinputsegmentremainingdistance}{0}}
		\pgfpathcircle{\pgfpoint{0.25\pgfdecoratedinputsegmentremainingdistance}{0}}%
		{0.05\pgfdecoratedinputsegmentremainingdistance}
		\pgfpathlineto{\pgfpoint{0.75\pgfdecoratedinputsegmentremainingdistance}{\pgfdecorationsegmentlength}}
		\pgfpathcircle{\pgfpoint{0.75\pgfdecoratedinputsegmentremainingdistance}{0}}%
		{0.05\pgfdecoratedinputsegmentremainingdistance}
		\pgfpathlineto{\pgfpoint{\pgfdecoratedinputsegmentremainingdistance}{0}}
		\pgfpathcircle{\pgfpoint{\pgfdecoratedinputsegmentremainingdistance}{0}}%
		{0.05\pgfdecoratedinputsegmentremainingdistance}
	}
	\state{final}
	{
		\pgfpathmoveto{\pgfpointdecoratedpathlast}
	}
}
\pgfmathdeclarefunction{gauss}{2}{%
	\pgfmathparse{1/(#2*sqrt(2*pi))*exp(-((x-#1)^2)/(2*#2^2))}%
}
\DeclareMathOperator*{\minimize}{minimize}
\DeclareMathOperator*{\sbjto}{subject\ to}
\DeclareMathOperator*{\blkdiag}{bdiag}

\DeclareMathOperator{\sat}{sat}
\DeclareMathOperator{\trace}{tr}

\renewcommand{\leq}{\leqslant}

\renewcommand{\geq}{\geqslant}

\newcommand{\R}{\mathds{R}}
\newcommand{\Nz}{\mathds{N}_0}
\newcommand{\N}{\mathds{Z}_{+}}
\newcommand{\EE}{\mathds{E}}
\newcommand{\PP}{\mathds{P}}

\newcommand{\bmat}[1]{\begin{bmatrix}#1\end{bmatrix}}
\newcommand{\abs}[1]{\left|#1\right|}
\newcommand{\norm}[1]{\left\|#1\right\|}
\newcommand{\secref}[1]{Section \ref{#1}}

\newcommand{\transp}{^\top}

\newcommand{\inverse}{^{-1}}

\newcommand{\zeros}{\mathbf{0}}

\newcommand{\st}{x}

\newcommand{\stest}{\tilde{x}}
\newcommand{\stfilt}{\hat{x}}
\newcommand{\estError}{\tilde{e}}
\newcommand{\dortho}{d_o}
\newcommand{\dschur}{d_s}
\newcommand{\stinit}{\st_0}

\newcommand{\A}{A}
\newcommand{\calA}{\mathcal{A}}
\newcommand{\Aortho}{\A_o}
\newcommand{\Aschur}{\A_s}
\newcommand{\B}{B}

\newcommand{\Bortho}{\B_o}
\newcommand{\Bschur}{\B_s}
\newcommand{\meas}{y}
\newcommand{\C}{C}
\newcommand{\control}{u}
\newcommand{\controlset}{\mathds{U}}
\newcommand{\wnoise}{w}

\newcommand{\cnoise}{\nu}
\newcommand{\snoise}{s}

\newcommand{\mnoise}{\varsigma}
\newcommand{\calB}{\mathcal{B}}
\newcommand{\calD}{\mathcal{D}}
\newcommand{\calC}{\mathcal{C}}
\newcommand{\calF}{\mathcal{F}}
\newcommand{\calG}{\mathcal{G}_t}
\newcommand{\calH}{\mathcal{H}}
\newcommand{\calQ}{\mathcal{Q}}
\newcommand{\calR}{\mathcal{R}}

\newcommand{\calS}{\mathcal{S}_t}

\newcommand{\calO}{\mathcal{O}}
\newcommand{\sigalg}{\mathfrak{F}}

\newcommand{\innovation}{\mathcal{I}}
\newcommand{\noisyInnovation}{\tilde{\mathcal{I}}}

\newcommand{\YY}{\mathfrak Y}
\newcommand{\XX}{\mathfrak Y^s}
\newcommand{\lra}{\longrightarrow}

\newcommand{\ee}{\mathfrak{\psi}}
\newcommand{\offset}{\boldsymbol{\eta}}
\newcommand{\gain}{\boldsymbol{\Theta}}
\newcommand{\authority}{u_{\max}}
\newcommand{\reachindex}{\kappa}
\newcommand{\reachab}{\mathrm{R}}
\newcommand{\Let}{\coloneqq}
\newcommand{\teL}{\eqqcolon}

\newtheorem{remark}{Remark}
\newtheorem{pstatement}{\sc Problem Statement}

\allowdisplaybreaks
\begin{document}

\begin{frontmatter}

\title{Stochastic Predictive Control under Intermittent Observations and Unreliable Actions} % Title, preferably not more 

\author[UIUC]{Prabhat K. Mishra} 
\author[IIT]{Debasish Chatterjee}
\author[Paderborn]{Daniel E. Quevedo}
\address[UIUC]{Coordinated Science Laboratory, University of Illinois at Urbana-Champaign, USA. \tt{m.prabhat}@outlook.com} 
\address[IIT]{Systems \& Control Engineering,
	Indian Institute of Technology Bombay, 
	India.  \tt{dchatter}@iitb.ac.in}
\address[Paderborn]{Department of Electrical Engineering (EIM-E), Paderborn University
	, Germany. \tt{dquevedo}@ieee.org}

                             % keyword list or with the 
                                          % help of the Automatica 
                                          % keyword wizard

\begin{abstract}                          % Abstract of not more than 200 words.
We propose a provably stabilizing and tractable approach for control of constrained linear systems under intermittent observations and unreliable transmissions of control commands. A smart sensor equipped with a Kalman filter is employed for the estimation of the states from incomplete and corrupt measurements, and an estimator at the controller side optimally feeds the intermittently received sensor data to the controller. The remote controller iteratively solves constrained stochastic optimal control problems and transmits the control commands according to a carefully designed transmission protocol through an unreliable channel. We present a (globally) recursively feasible quadratic program, which is solved online to yield a stabilizing controller for Lyapunov stable linear time invariant systems under any positive bound on control values and any non-zero transmission probabilities of Bernoulli channels.    
\end{abstract}

\end{frontmatter}

\section{Introduction}
Predictive techniques for networked control systems (NCSs) have significantly advanced over the past decade. Emphasis has been placed on information loss \cite{varutti2009compensating, grune2009prediction}, recursive feasibility issue \cite{pin2011network}, stability analysis  \cite{Quevedo2011, Quevedo-12}, and focus on applications such as drinking water networks \cite{pereira2016application, sampathirao2017gpu}. In addition, NCSs have been constructed in diverse, interesting and important applications including haptic collaboration over the internet \cite{shirmohammadi2004evaluating, hespanha2000haptic, hikichi2002evaluation}, building automation \cite{newman1996integrating}, vehicle control \cite{seiler2001analysis}, mobile sensor networks \cite{ogren2004cooperative} to name only a few. Predictive techniques make the controller capable of handling constraints while optimizing a desired performance objective.  
\begin{figure*}
	\centering
	\begin{adjustbox}{width = 0.8\textwidth}
		%\begin{figure*}
%	\begin{adjustbox}{width = \textwidth}
		%		\input{blockdia}
		\begin{tikzpicture}
		\tikzstyle{pinstyle} = [pin edge={to-,thin,black}]	
		\tikzstyle{block} = [draw, fill=blue!10, rectangle, 
		minimum height=2em, minimum width=1]
		\tikzstyle{blockgreenl} = [draw, fill=green!10, rectangle, 
		minimum height=2em, minimum width=1cm]
		\tikzstyle{blockgreen} = [draw, fill=green!10, rectangle, 
		minimum height=2em, minimum width=1cm] 
		\tikzstyle{blockred} = [draw, fill=red!10, rectangle, 
		minimum height=2em, minimum width=1.1cm] 
		\tikzstyle{blockcover} = [draw, fill=white, rectangle, dashed,  
		minimum height=2em, minimum width=1]                        
		%\tikzset{cross/.style={cross out, draw=black, minimum size=2*(#1-\pgflinewidth), inner sep=0pt, outer sep=0pt},					
		\tikzstyle{sum} = [draw, circle]					     
		%\tikzstyle{pinstyle} = [pin edge={to-,thin,black}]	
		%===================================================
		\node[coordinate] (0) at (0,0) {};
		\node [starburst, right= 0 cm of 0 , fill=white,draw, dashed, pin={[pinstyle]above:$\cnoise_t$}] (ControlChannel) {Control channel};
		%\node [blockred, right= 2cm of 0, pin={[pinstyle]above:$\wnoise_t$}] (plant) {$\st_{t+1}=A\st_t + B \control_t$};	
		\node [blockred, right= 1cm of ControlChannel] (actuator) {Actuator};	
		\node [blockred, right= 1cm of actuator, pin={[pinstyle]above:$\wnoise_t$}] (plant) {Plant};		
		\node [blockgreen, below= 2cm of ControlChannel] (controller) {Controller};
		\node [blockgreen, right= 2cm of controller] (estimator) {Estimator};
		\node [sum, right= 1cm of plant, pin={[pinstyle]above:$\mnoise_{t+1}$}] (sensor) {+};
		%\node [blockred, right= 1cm of plant, pin={[pinstyle]above:$\mnoise_t$}] (sensor) {Sensor};
		\node [blockred, right= 1cm of sensor] (filter) {Filter};
		\node [starburst, below = 1.65 cm of filter, fill=white,draw, dashed, pin={[pinstyle]below:$\snoise_t$}] (SensorChannel) {Sensor channel};
		\draw[very thin ,->] (estimator.west) to node[auto, swap] {$\stest_{t}, \noisyInnovation_t$}  (controller.east);
		\node [draw, dotted, fit= (controller) (estimator), inner sep = 5pt] (G) {};
		\node [draw, dotted, fit= (actuator) (plant) (sensor) (filter), inner sep = 5pt] (H) {};
		%\node[anchor=north, right= 0.05cm of plant.20, blue] {$C \st_{t+1}$};
		
		%\node[anchor=west, below= 1cm of controller.south, blue] {$\hat{\st}_t$,$\hat{\meas}_t$,$\meas_t$};
		%\node[anchor=west, right= 0.4cm of actuator.-15, blue] {$\control_t^a$};
		%\node[anchor=east, right= 0.2cm of sensor.30, blue] {$\meas_{t+1}$};
		\node[anchor=east, above= 0.5cm of SensorChannel.north, xshift = 0.5cm,  blue] {$\stfilt_{t}, \meas_t$};
		%\node[anchor=west, right= 0.1cm of controller.15, blue] {$\stest_{t}, \noisyInnovation_t$};
%		\node[anchor=east, below = 0.6cm of filter.5, blue] {$\stfilt_{t+1}, \meas_{t+1}$};
		% \node[anchor=north, above = 0.4cm of ($(actuator.north west)!.8!(filter.north east)$), blue] {$\control_{t}$};
		%\node[anchor=west, above= 0.8cm of filter.west, blue] {$\control_t$};
		%\node[anchor=north, right= 0.2cm of plant.20, blue] {$C \st_t$};
		%==========================================
		%\draw[very thin ,->] (actuator.east) edge[bend right] (filter.west);
		%\draw[very thin ,->] (actuator.north) to[myncbar,angle=90,arm=1.5] (filter.north);
		\draw[very thin ,->] (actuator.north) -- ($(actuator.north) + (0,1.5)$) to node[auto, swap] {$\control_t^a$} ($(filter.north) + (0,1.5)$) -- (filter.north);
		\draw[very thin ,->] (plant.east) to node[auto, swap] {$C \st_{t+1}$}  (sensor.west);
		\draw[very thin ,->] (sensor.east) to node[auto, swap] {$\meas_{t+1}$}  (filter.west);
		\draw[very thin ,->] (filter.south) -- (SensorChannel.north);
		\draw[very thin ,->] (ControlChannel.east) -- (actuator.west);
		\draw[very thin ,->] (actuator.east) to node[auto, swap] {$\control_t^a$} (plant.west);
		\draw[very thin ,->] (SensorChannel.west) -- (estimator.east);		
		\draw[very thin ,->] (controller.north) -- (ControlChannel.south);
		\end{tikzpicture}
%	\end{adjustbox}	
%	\caption[Block diagram]{The red blocks are situated at the plant and the green blocks at the remote controller. The red and the green blocks communicate via erasure channels. The Kalman filter employs the measurement $\meas_t$ and the applied control $\control_t^a$ to estimate the state $\stfilt_{t+1}$. The acknowledgements of successful transmissions of control commands are causally available to the controller and the estimator. The estimator is also aware of the sensor channel dropout $\snoise_t$ at time $t$.}		
%%	\label{Fig:Setup}		
%\end{figure*}
	\end{adjustbox}	
	\caption{The red blocks are situated at the plant and the green blocks at the remote controller. The red and the green blocks communicate via erasure channels. The Kalman filter employs the measurement $\meas_t$ and the applied control $\control_t^a$ to estimate the state $\stfilt_{t+1}$. The acknowledgements of successful transmissions of control commands are causally available to the controller and the estimator. The estimator is also aware of the sensor channel dropout $\snoise_t$ at time $t$.}		
	\label{Fig:Setup}		
\end{figure*}
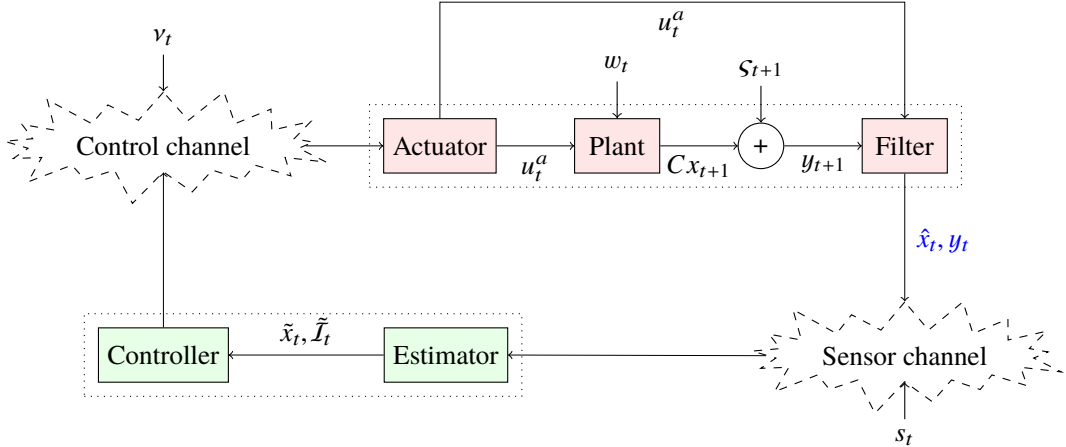
\par Most of the predictive control techniques neglect uncertainties when designing the objective function,  e.g., \cite{Quevedo2011}, or consider only the worst case scenario over uncertainties, e.g., \cite{pereira2016application}. In both cases the resulting controllers are conservative and do not take advantage of the available statistics. However, while controller design without incorporating such statistical information is easier and simpler than otherwise, this results in a decline in the desired performance \cite{ref:PDQ-15}. Stochastic predictive controllers offer a way out of such conservatism. Two major challenges need to be overcome in the stochastic controller design -- tractability of the underlying constrained stochastic optimal control problem (CSOCP) and guaranteeing stability in some suitable sense. We refer the readers to \cite{OutBitmead, heirung2018stochastic} for discussions on the underlying challenges in stochastic predictive control under the settings of perfect channels. Recently, \cite{ref:PDQ-15, ref:amin-10, hokayem12network} designed tractable and stabilizing controllers for networked systems that take the probability distributions of the uncertainties and network induced effects into account at the synthesis stage. However, \cite{ref:PDQ-15, ref:amin-10, hokayem12network} rely on the assumption that the sensor channel is perfect and there are dropouts only in the control channel. This assumption is reasonable for applications as mentioned in \cite{QuevedoMultiple16,wuStability-10}, but in a wide range of applications available sensor-communication channels are unreliable. Therefore, an extension of \cite{ref:PDQ-15, ref:amin-10, hokayem12network} to the setting of an unreliable sensor channel is important for a large class of applications, and lies at the heart of the current work.   
\par Apart from unreliable communication channels, typically measurements are corrupted by sensor noise and/or full state information are not available. Moreover, since actuators are physical devices, constraints on the control actions must be satisfied for all realizations of the uncertainties.  
Although the synthesis of stabilizing constrained control under incomplete and corrupt observations is an interesting problem, the literature on this topic is sparse, apart from \cite{ref:Hokayem-12, PDQ-LCSS}. In both of these articles a stochastic predictive controller was proposed over ideal communication channels with an affine saturated innovation feedback policy employing a Kalman filter and constant negative drift conditions.
However, their implementation over unreliable networks is non-trivial due to boundedness issues in Kalman filtering with intermittent observations \cite{intermittent, quevedo2013TAC03}. 

\par Several interesting methods have been proposed to estimate the system states over unreliable sensor channels \cite{bemporad2010networked, intermittent, quevedo2013TAC03, sun2008optimal, liu2013moving, jin2006state, gupta2009data, plarre2009kalman, kar2012kalman, Quevedo17TAC, Murray10TAC, gupta2007optimal}. However, to the best of our knowledge the case of sensor channel erasure under the settings of stochastic predictive constrained control has not been examined so far due to, as already mentioned, boundedness issues involved in Kalman filtering. To be precise, the conditional error covariance matrix exhibits unbounded oscillations almost surely, and is therefore computationally difficult to deal with. 
In this article we adapt the idea of smart sensors \cite{Murray10TAC, gupta2007optimal} and present a controller that has a computationally tractable underlying optimal control problem. The proposed controller  ensures stability of the closed-loop states in a suitable sense under any positive bound on control and any successful transmission probabilities of the sensor and the control channels. Our results hold for the largest class (to date) of linear time invariant (LTI) discrete dynamical systems known to be stabilizable under bounded control actions even without dropouts. For that purpose, we generalize the approach of \cite{ref:PDQ-15, ref:amin-10, hokayem12network, ref:Hokayem-12, PDQ-LCSS} by considering the unreliable sensor channel in the network system architecture. 
\par The main features and contributions of the proposed approach are as follows:
\begin{itemize}[leftmargin = *]
	\item We consider an LTI system, with incomplete and corrupt measurements, remotely controlled over unreliable channels.
	\item We present a tractable and recursively feasible quadratic program to be solved periodically online, which provides optimal control that minimizes the expected quadratic cost function and ensures a good closed-loop behaviour of states.
	\item In order to compensate the effect of packet dropouts in the sensor channel, we employ a remote estimator. We show that the remote estimator provides the conditional expectation of states given causally available information at the estimator.
	\item We employ a novel class of feedback policies and show convexity of the underlying optimization program.  
	\item For systems having all eigenvalues inside the unit circle and semi-simple eigenvalues on the boundary (if any), we show mean-square boundedness of the controlled system for any positive bound on the control actions. Mean-square boundedness of the controlled states is guaranteed for any non-zero probabilities of successful transmissions in both the sensor channel and the control channel.
\end{itemize}      
\par This article exposes as follows: In \secref{s:Problem Setup} we discuss the basic assumptions and formally define the problem statement. The system setup is presented in \secref{s:setup}. We develop our main results on tractability and stability in \secref{s:tractability} and \secref{s:stability}, respectively. We validate our theoretical results by numerical experiments in \secref{s:experiment} and conclude in \secref{s:epilogue}. Some proofs are presented in the Appendix.
\subsection*{Notation}
Let $\R, \Nz, \N$ denote the set of real numbers, the non-negative integers and the positive integers, respectively. We use the symbol \( \zeros \) to denote a matrix of appropriate dimensions with all elements $0$.  
For any vector sequence \((v_n)_{n\in\Nz}\) taking values in some Euclidean space, let \( v_{n:k} \) denote the vector \(\bmat{v_n\transp & v_{n+1}\transp & \cdots & v_{n+k-1}\transp}\transp\), \(k\in\N\). The notations \(\EE_z[\cdot]\) and \(\EE[\cdot \mid z ]\) are interchangeably used for the conditional expectation with given \(z\). For a vector \(v\), its $i^{th}$ element is denoted by $v^{(i)}$. 
Similarly, $M^{(j,:)}$ denotes the $j^{th}$ row of a given matrix $M$. 
Let $\sigma_1(M)$ denote the largest singular value of $M$, and $M^\dagger$ its Moore-Penrose pseudo inverse. A block diagonal matrix $M$ with diagonal entries $M_1, \cdots , M_n$ is represented as $M = \blkdiag\{M_1, \cdots , M_n\}$ and \(I_d\) is the \(d\times d\) identity matrix. For a real quantity $\xi$, its positive component $\xi_+$ and negative component $\xi_-$ are defined to be $ \max\{0,\xi\}$ and $ \max\{0,-\xi\}$, respectively. For a given positive semi-definite (or definite) matrix $P$ and a vector $v$, the notation $\norm{v}_P^2$ is used to denote the scalar $v\transp P v$.

\section{Problem Statement}\label{s:Problem Setup}
Let us consider a discrete time dynamical system 
\begin{subequations}\label{e:system}
	\begin{align}
	\st_{t+1} &= \A\st_t + \B\control_t^a + \wnoise_t  \label{e:steq} \\ 
	\meas_t & = \C \st_t + \mnoise_t,
	\end{align}
\end{subequations}
where $t \in \Nz$, and $\st_t \in \R^d$, $\control_t^a \in \R^m$, $\meas_t \in \R^q$ are the states, the applied control to the plant and the measurements, respectively, at time $t$. The additive process noise $\wnoise_t \in \R^d$ and the measurement noise $\mnoise_t \in \R^q$ are zero-mean Gaussian. System matrices $\A,\B$ and $\C$ are known matrices of appropriate dimensions and the matrix pair $(\A,\B)$ is stabilizable.  
At each time $t$ the control $\control_t^a$ is constrained to take values in the admissible set
\begin{equation}\label{e:controlset}
\controlset \Let \{v\in\R^m\mid \norm{v}_\infty \leq \authority \}, 
\end{equation}		
where the (uniform) bound $\control_{\max} > 0$ is preassigned.

\begin{remark}
\rm{
Admissible control set of the form
\begin{equation*}\label{e:ControlSet}
\controlset^{\prime} \Let \Biggl\{ v \in \R^m \Biggm| \abs{v^{(i)}} \leq U_i \text{ for } i =1, \ldots, m \Biggr\},
\end{equation*}  
for not necessarily equal values $U_i$, can be transformed easily into $\controlset$ as in \eqref{e:controlset}. Please see \cite[Remark 1]{PDQ_Policy} for more details.
}
\end{remark}
The sensor channel and the control channel both are unreliable with successful transmission probabilities, $p_s$ and $p_c$, respectively. We represent the dropout at time $t$ in the sensor and the control channel by i.i.d. Bernoulli random variables $\snoise_t$ and $\cnoise_t$, respectively, see Fig.\ \ref{Fig:Setup}.
	We have the following assumption:
	\begin{enumerate}[label={\rm (A\arabic*)}, leftmargin= *, widest=3, align=left, start=1, nosep]
		\item \label{as:uncorrelated} The dropout processes $(\cnoise_t)_{t \in \Nz}$ and $(\snoise_t)_{t \in \Nz}$ are mutually independent and individually i.\ i.\ d. They are also independent of the process noise and the measurement noise processes $(\wnoise_t)_{t \in \Nz}$ and $(\mnoise_t)_{t \in \Nz}$, respectively. 
	\end{enumerate}
The inclusion of the constraint \eqref{e:controlset} at the synthesis stage is achieved by the following constrained stochastic optimal control problem (CSOCP) that is solved iteratively over time, and that constitutes the backbone of predictive control techniques:
\begin{equation}
\label{e:opt control problem}
\begin{aligned}
& \minimize_{\text{control policies}}	&&  \text{ a quadratic objective function }\\
& \sbjto	 & &  \begin{cases}				
\text{system dynamics }\eqref{e:system},\\
\text{hard constraint on control } \eqref{e:controlset}.
\end{cases}
\end{aligned}
\end{equation}			
It is well known that the stochastic optimal control problem in the form of \eqref{e:opt control problem} is computationally intractable even in the presence of the full state information and perfect channels. An affine feedback policy based approach is often used to present a tractable surrogate of CSOCP \eqref{e:opt control problem} when an expected quadratic cost is used as objective function \cite{mesbah_16_survey}. 
In order to mitigate the effect of incomplete and corrupt measurements, a Kalman filter (see \secref{s:KF}) at the sensor with the following assumptions is employed.				
	\begin{enumerate}[label={\rm (A\arabic*)}, leftmargin= *, widest=3, align=left, start=2, nosep]
		\rm
		\item The matrix pair $(\A,\C)$ is observable. \label{as:observable}
		\item The initial condition $\stinit$, the process and the measurement noise vectors are normally distributed and mutually independent, i.e. $\st_0 \sim N(0,\Sigma_{\st_0})$, $\wnoise_t \sim N(0,\Sigma_{\wnoise}), \mnoise_t \sim N(0,\Sigma_{\mnoise})$, with $\Sigma_{\st_0} \succeq 0$, $ \Sigma_{\wnoise} \succeq 0$ and $\Sigma_{\mnoise} \succ 0$. \label{as:normal distribution} 
		\item The matrix pair $(\A, \Sigma_{\wnoise}^{1/2})$ is controllable.\label{as:stationaryP}
	\end{enumerate}
At the controller end, an optimal estimator is employed to mitigate the effects of sensor channel dropouts; see Fig. \ref{Fig:Setup} for a schematic. In one hop communication links, it is a standard practice to transmit error free receipt acknowledgements (or negative acknowledgements) \cite[page 207]{kuroseNetworking}. Therefore, the estimator is aware of the previously applied controls to the plant. Accordingly, we have the following assumption: 
%\begin{assumptionnn}\mbox{}
	\begin{enumerate}[label={\rm (A\arabic*)}, leftmargin= *, widest=3, align=left, start=5, nosep]
		\item The acknowledgements of the successfully transmitted control commands are causally available to the controller and the estimator. \label{as:tcp} 		
	\end{enumerate}	
%\end{assumptionnn}
The objective of this paper is to present a tractable surrogate of CSOCP \eqref{e:opt control problem} and ensure stability (in some suitable sense) of the closed-loop system \eqref{e:steq} when control is computed by iteratively solving \eqref{e:opt control problem}. Since in the examined situation asymptotic stability of the origin cannot be achieved due to the unbounded support of the additive process noise and the measurement error, we are interested in the notion of mean square boundedness, as defined below:
\begin{defn}[Mean square boundedness]
	An $\R^d$-valued random process $(\st_t)_{t \in \Nz}$ is said to be mean-square bounded (MSB) if there exists some $\gamma > 0$ such that \[ \sup_{t \in \Nz}\EE_{\XX_0}[\norm{\st_t}^2] \leq \gamma,\]
	where $\XX_0$ is the information available at $t=0$. 
\end{defn}  
It is well known that, even with perfect communication channels, an LTI system of the form \eqref{e:system} with system matrix $\A$ having eigenvalues outside the unit circle cannot be (globally) stabilized by any control technique with help of bounded control actions \cite[Abstract]{Sussmann97}, \cite[Theorem 1.7]{ref:ChaRamHokLyg-12}. Consequently, for the present networked case, we impose the following assumption:
	\begin{enumerate}[label={\rm (A\arabic*)}, leftmargin= *, widest=3, align=left, start=6, nosep]
		\item \label{as:lyapunov} The system matrix $\A$ has all eigenvalues on the unit disk and those on the unit circle are semi-simple.\footnote{The stabilizability of LTI systems with non-semi-simple eigenvalues on the unit circle under bounded control actions remains an open problem \cite{ref:ChaRamHokLyg-12}. Therefore, the considered class of systems is the largest class of LTI systems, known till date, stabilizable under bounded control actions.}   		
	\end{enumerate}		       
Dynamical systems satisfying Assumption \ref{as:lyapunov} are well-studied as Lyapunov stable systems \cite{ref:ChaHokLyg-11}.
In the presence of unbounded disturbances, an open-loop Lyapunov stable (but not asymptotically stable) system leads to unstable trajectories with probability one, unless the control is designed with sufficient care. Such designs are non-trivial if hard constraints on the control have to be satisfied at all times. Traditional MPC for LTI systems under hard bounds on control inputs has been widely studied \cite{borrelli2011predictive}. Unfortunately, the approaches in \cite{borrelli2011predictive} do not guarantee mean square boundedness in the presence of unbounded disturbances. The control strategy developed in the present article is for stochastic systems (additive disturbance has unbounded support as well as observations are intermittent and control commands are unreliable) where most of the well developed tools of deterministic MPC do not carry over. In particular, tools developed assuming deterministic settings mostly rely on terminal set and cost methods \cite{ref:rawlings-09, Quevedo2011}. Further, the construction of suitable positively invariant sets in the presence of disturbances with unbounded support is impossible \cite[\S3]{ref:May-14}; please also see \cite[Lemma 1]{huang-bitmead15}. In order to transcend beyond the regime of terminal set and cost method drift conditions are used for stochastic systems. We do not invoke martingale arguments directly for our drift conditions. Instead, our result is based on the approach in \cite[Theorem 1]{ref:PemRos-99}. The latter work contains a delicate proof of their main result relying on the Burkholder's inequality, Doob decomposition in the theory of martingales and some other arguments. We recall the following result:
\begin{thm}[{\cite[Theorem 1, Corollary 2]{ref:PemRos-99}}]
	\label{t:PemRos-99}
	Let \((X_t)_{t\in\Nz}\) be a family of real valued random variables on a probability space \((\Omega, \sigalg, \PP)\), adapted to a filtration \((\sigalg_t)_{t\in\Nz}\). Suppose that there exist scalars \(a, b, c > 0\) such that 
	\begin{equation*}
	\begin{aligned}
	& \EE_{\sigalg_t}[X_{t+1} - X_t]  \leq -a \quad\text{on the event } X_t > b, \\
	& \EE\bigl[\abs{X_{t+1} - X_t}^4\,\big|\, X_0, \ldots, X_t\bigr]  \leq c\quad\text{for all }t\in\Nz.
	\end{aligned}
	\end{equation*}			
	Then there exists a constant \(\gamma > 0\) such that \[\sup_{t\in\Nz}\EE\bigl[ \left((X_t)_+\right)^2 \mid  \sigalg_0 \bigr] \leq \gamma .\]
\end{thm}
The first condition of the above theorem is called the \emph{constant negative drift condition} and it is active when $X_t$ is larger than some $b > 0$. The second condition is called \emph{skip-free condition} and is needed to avoid long jumps when $X_t \leq b$. 
The above theorem gives sufficient conditions for the mean square boundedness of the positive component $(X_t)_+$ of a scalar process $(X_t)_{t \in \Nz}$ and extended for the vector processes in \cite{ref:ChaRamHokLyg-12} utilizing the Assumption \ref{as:lyapunov} and decomposition of the system dynamics into orthogonal and Schur stable subsystems, which is given explicitly in \secref{s:stability}. We utilized the above theorem and the idea of decomposition in our previous works \cite{ref:PDQ-15, PDQ_Policy} in a limited context. A part of our stability result is along the lines of \cite[Theorem 1.2]{ref:ChaRamHokLyg-12} but satisfaction of the second condition of Theorem \ref{t:PemRos-99} is non-trivial in the setting of the present article. In particular, \cite[Theorem 1.2]{ref:ChaRamHokLyg-12} shows the existence of a stabilizing history dependent feedback policy by considering a $\reachindex-$subsampled process when the system matrix $A$ is orthogonal with reachability index $\reachindex$ and channels are perfect. In this article we extend the stability analysis of \cite{ref:ChaRamHokLyg-12,ref:RamChaMilHokLyg-10} under the settings of unreliable channels.   
\begin{remark}
	\rm{
		The	recursive feasibility of stochastic predictive control techniques under state constraints is challenging whenever involved noise processes have unbounded support \cite{primbs2009stochastic}. The inclusion of state constraints within our framework can be investigated along the lines of \cite{hokayem10OutState, RE-SPC}. For simplicity of the presentation we have not considered state constraints in this article.		
	}
\end{remark}
The problem statement of the present article is formally given below:\\ 
\begin{pstatement}
	\label{ps:main statement}
	\rm{
		Present a tractable, stabilizing and recursively feasible surrogate of CSOCP \eqref{e:opt control problem} under the assumptions \ref{as:uncorrelated} -- \ref{as:lyapunov}.	}  
\end{pstatement}
\section{Setup}\label{s:setup}
As illustrated in Fig. \ref{Fig:Setup}, we employ a Kalman filter at the sensor and an estimator at the controller.	The information of the past outputs and previously applied control is available at the filter. At each time $t$ the filtered state $\stfilt_t$ and output $\meas_t$ are transmitted through the sensor channel. Since, the sensor channel is affected by Bernoulli dropouts, either the transmitted information reaches the estimator or it is lost. For the sake of computational tractability we consider quadratic cost functions, which are minimized over a class of policies. Control commands obtained by solving the optimization programs are transmitted through an erasure channel. To mitigate the effects of dropouts in the control channel we employ ad-hoc  transmission strategies. A detailed discussion on the above features is presented below:  

\subsection{Expected quadratic cost}
Let us fix an optimization horizon $N \in \N$ and recalculation interval (control horizon) $N_r \leq N$. Let $\XX_t$ be the information available at the controller/estimator at the time of optimization $t$ (see \secref{s:optimal estimator} for a precise definition). The cost $V_t$ is defined to be the conditional expectation of the quadratic cost in one optimization horizon with given information $\XX_t$. We define $V_t$ as 
\begin{equation}\label{e:cost}
V_t \Let \EE_{\XX_t} \left[ \sum_{k = 0}^{N-1} (\norm{\st_{t+k}}^2_{Q} + \norm{\control_{t+k}^a}^2_{R}  ) + \norm{\st_{t+N}}^2_{Q_N}\right],
\end{equation}
where, $Q, Q_N$ are given symmetric positive semi-definite matrices of appropriate dimensions and $R$ is given symmetric positive definite matrix.	
The compact form representation of the system \eqref{e:system} over one optimization horizon is as follows:
\begin{subequations}\label{e:stacked dynamics}
	\begin{align}
	\st_{t:N+1} &= \calA\st_t + \calB\control_{t:N}^a + \calD \wnoise_{t:N} \label{e:stacked state}\\
	\meas_{t:N+1} &= \calC \st_{t:N+1} + \mnoise_{t:N+1},
	\end{align}
\end{subequations}
where $\calA$, $\calB$, $\calC$ and $\calD$ are standard matrices of appropriate dimensions. 
The cost function \eqref{e:cost} can also be written in a compact form for later use as follows:
\begin{equation}\label{e:cost_compact}
V_t = \EE_{\XX_t} \left[ \norm{\st_{t:N+1}}^2_{\calQ} + \norm{\control_{t:N}^a}^2_{\calR} \right],
\end{equation}	
where $\calQ$ and $\calR$ are standard block diagonal matrices of appropriate dimensions.

\subsection{Kalman filter}\label{s:KF}
In this section we recall the framework of stochastic predictive control using Kalman filtering. The detailed discussion is available in \cite{ref:Hokayem-12}. For each $t$ let 
$ \YY_t \Let \{ \meas_0, \cdots, \meas_t, \control_0^a, \cdots, \control_{t-1}^a \} $
denote the set of observations up to time $t$.
For $t,s \in \Nz$, $t \geq s$, let us define $\stfilt_{t \mid s} \Let \EE_{\YY_s}\left[ \st_t \right]$ and $P_{t \mid s} \Let \EE_{\YY_s} \left[ (\st_t - \stfilt_{t \mid s}) (\st_t - \stfilt_{t \mid s}) \transp \right]$, and for brevity of notation, we denote $\stfilt_{t \mid t}$ by $\stfilt_t$ and $P_{t \mid t}$ by $P_t$. We need the following result related to Kalman filtering for which recursions are \cite[p.102]{kumar1986stochastic}:
\begin{equation}\label{e:filter dynamics}
\begin{aligned}
\stfilt_{t+1} &=  \stfilt_{t+1 \mid t} + K_t (\meas_{t+1}-\C\stfilt_{t+1 \mid t}) \\
P_{t+1} &= P_{t+1 \mid t} - K_t \C P_{t+1 \mid t}
\end{aligned}
\end{equation}
where $\stfilt_{t+1 \mid t} = \A \stfilt_t + \B \control_t^a$, $P_{t+1 \mid t} = \A P_t \A \transp + \Sigma_{\wnoise} $ and $ K_t = P_{t+1 \mid t}\C \transp \Bigl(\C P_{t+1 \mid t}\C\transp + \Sigma_{\mnoise} \Bigr)\inverse$.
We intialize the Kalman filter by setting $\stfilt_{0 \mid -1} = \zeros, P_{0 \mid -1} = \Sigma_{\st_0} $, we get $\stfilt_0 = K_0(\C\st_0 + \mnoise_0)$. Let us recall the results of \cite[Lemma 8]{ref:Hokayem-12} and \cite[Lemma 4.2.2]{balakrishnan1987kalman} that there exists a constant $\rho > 0 $ such that 
\begin{equation}
\EE_{\YY_t} \left[ \norm{\st_t - \stfilt_t }^2 \right] \leq \rho \quad \text{ for all } t.
\end{equation} 
Let us define filtered output $\hat{\meas}_t \Let \C \stfilt_t$, the innovation term $\innovation_t \Let \meas_t - \hat{\meas}_t$ and Kalman filter error $e_t \Let \st_t - \stfilt_t$. 
A straight-forward calculation gives the innovation term for one optimization horizon as follows:
\begin{equation} \label{e:est error}
\innovation_{t:N+1} = \calC\calF_t e_t + \calO_t \wnoise_{t:N} + (I-\calC \calH_t) \mnoise_{t:N+1}
\end{equation}
where 
\begin{align*}
\calF_t & \Let \bmat{I_d \\ \phi_t \\ \phi_{t+1}\phi_t  \\ \vdots \\ \phi_{t+N-1}\cdots \phi_t } \\
\calO_t & \Let \bmat{0 & \cdots & 0 & 0 \\ \Gamma_t & \cdots & 0 & 0 \\ \phi_{t+1}\Gamma_t & \cdots & 0 & 0 \\ \vdots & \cdots & \vdots & \vdots \\ \phi_{t+N-2} \cdots \phi_{t+1}\Gamma_t & \cdots & \Gamma_{t+N-1} & 0 \\ \phi_{t+N-1}\cdots \phi_{t+1}\Gamma_t & \cdots & \phi_{t+N-1}\Gamma_{t+N-2} & \Gamma_{t+N-1} } \\
\calH_t & \Let \bmat{0 & 0 & \cdots & 0 & 0 \\0 & K_t & \cdots & 0 & 0 \\ 0 & \phi_{t+1} K_t & \cdots & 0 & 0 \\ 0 & \vdots & \cdots & \vdots & \vdots \\ 0 & \phi_{t+N-2} \cdots \phi_{t+1} K_t & \cdots & K_{t+N-1} & 0  \\ 0 & \phi_{t+N-1}\cdots \phi_{t+1} K_t & \cdots & \phi_{t+N-1} K_{t+N-2} & K_{t+N-1} } 
\end{align*}
The dynamics of $\stfilt_t$ can be given by
\begin{equation}\label{e:estimator st}
\stfilt_{t+1} = \A \stfilt_t + \B \control_t^a + \hat{\wnoise}_t,
\end{equation}
where $\hat{\wnoise}_t = K_t(C\A e_t + C \wnoise_t + \mnoise_{t+1})$. 
The above equations \eqref{e:est error} and \eqref{e:estimator st} are standard in literature and can be found by using \eqref{e:system} and \eqref{e:filter dynamics}. Since $\hat{\wnoise}_t$ is a linear sum of three mutually independent Gaussian random variables, it is also Gaussian with a time varying and bounded variance due to $K_t$. 
\subsection{Optimal estimator}\label{s:optimal estimator}
For each $t$ let $\mathcal{M}_t \Let \bigl\{ \snoise_t\meas_t, \snoise_t \stfilt_t, \snoise_t, \control_{t-1}^a \bigr\}$ 
and $ \XX_t \Let \bigl\{ \mathcal{M}_0, \ldots , \mathcal{M}_t  \bigr\} $ denote the set of data available at the estimator up to time $t$, with the convention that $\control_{-1}^a = \zeros$. 
We have the following results: 
\begin{lem}\label{lem:ApproxEstimate}
	Let $\stest_t \Let \EE \left[ \st_t \mid \XX_t  \right]$. Then  
	\begin{equation}\label{e:est}
	\stest_t = \snoise_t \stfilt_t + (1-\snoise_t) \left( \A \stest_{t-1} + \B \control_{t-1}^a \right).
	\end{equation}
\end{lem}
A proof of Lemma \ref{lem:ApproxEstimate} is given in the appendix. The estimator \eqref{e:est} has been widely used in the literature. The above Lemma \ref{lem:ApproxEstimate} proves that it is optimal under the settings of the present article.
We define the estimation error $\estError_t \Let \stfilt_t - \stest_t$ and initialize $\stest_{-1} = \zeros$, $\estError_{-1} = \zeros$. 
	
\subsection{Affine saturated received innovation feedback policy}
Affine feedback parametrizations in terms of the innovation sequence are standard in the literature \cite{bosgra2003}. Since the optimization is carried over a particular class of policies, the obtained solution is sub-optimal. In order
to satisfy hard bounds on the control actions while retaining computational tractability, affine saturated innovation feedback policies are used in \cite{ref:Hokayem-12, PDQ-LCSS}. Since innovation terms are affected by sensor channel dropouts under the settings of this article, a received innovation term is employed in the feedback parametrization used in the present article. We consider the following causal feedback policy class for $\ell = 0,\cdots , N-1$, 
\begin{equation}\label{e:policyclass}
\begin{aligned}
\control_{t+\ell} &= \eta_{t+\ell} + \sum_{i=0}^{\ell} \theta_{\ell,t+i}\ee_i (\snoise_{t+i}\innovation_{t+i}) \\
&= \eta_{t+\ell} + \sum_{i=0}^{\ell} \snoise_{t+i} \theta_{\ell,t+i}\ee_i (\innovation_{t+i}),
\end{aligned}
\end{equation}
where \(\ee_i:\R^q\lra\ \R^q\) is a measurable map for each \(i\) such that  $\norm{\ee_i(\meas_{t+i}-\hat{\meas}_{t+i})}_{\infty} \leq \ee_{\max}$.  
Let us denote $\noisyInnovation_t \Let \snoise_t \innovation_t$ for brevity.
The above control policy class \eqref{e:policyclass} can be represented in a compact form as follows:
\begin{equation} \label{e:policy out channel}
\control_{t:N} = \offset_t + \gain_t \ee(\noisyInnovation_{t:N})
\end{equation} 
where \(\offset_t\in\R^{m N}\), 
$\ee \Let \bmat{\ee_1\transp & \cdots & \ee_{N}\transp}$ and \(\gain_t\) is the following block triangular matrix
	\begin{equation} \label{e:gain}
	\gain_t = \bmat{ \theta_{0, t} & \zeros & \cdots & \zeros & \zeros \\ \theta_{1, t} & \theta_{1, t+1}  & \cdots & \zeros & \zeros \\ \vdots & \vdots & \vdots & \vdots & \vdots\\ \theta_{N-1, t} & \theta_{N-1, t+1} & \cdots & \theta_{N-1, t+N-2} & \theta_{N-1, t+N-1} },
	\end{equation}
	with each \(\theta_{k, \ell} \in \R^{m\times q}\) and  \[ \norm{\ee(\noisyInnovation_{t:N})}_{\infty} \leq \ee_{\max}. \] 
We choose $\ee_i$'s to be component-wise odd functions, e.g., standard saturation function, sigmoidal function, etc.

\begin{lem}\label{lem:hardControl}
	The hard constraint on control \eqref{e:controlset} under the class of control policies \eqref{e:policyclass} is equivalent to the following constraint:
	\begin{equation}\label{e:decisionboundsingle}
	\abs{\offset_t^{(i)}} +  \norm{\gain_t^{(i,:)}}_1\varphi_{\max}  \leq \authority
	\end{equation}		
\end{lem}
\begin{pf}
	In view of the dropout process $\cnoise_t \in \{0,1 \}$,
	\begin{align*}
	\norm{\control_{t:N}^a}_{\infty} \leq \authority & \iff \norm{\control_{t:N}}_{\infty} \leq \authority \\
	& \iff \norm{\offset_t + \gain_t \psi (\innovation_{t:N})}_{\infty} \leq \authority.
	\end{align*}
	Now the assertion follows from \cite[Proposition 3]{hokayem2009stochastic}. 
\end{pf}

\subsection{Transmission protocol}

In order to mitigate the effects of packet dropouts in the control channel, several transmission protocols are discussed in \cite{ref:PDQ-15}. We consider one of them, which is formally defined below. Our approach remains valid for the other protocols as well, provided minor adjustments are made.   
\begin{enumerate}[label={\rm (TP)}, leftmargin=*, widest=3, align=left, start=2]
	\item At the beginning of each optimization instant ($t = 0, N_r, 2N_r, \ldots $), the buffer is emptied. For $\ell \leq N_r-1 $, $\control_{t+\ell}$ is transmitted and directly applied to the plant if successfully received at the actuator. In additon, $\eta_{t+\ell +1}, \ldots , \eta_{t+N_r-1}$ are also transmitted until the first successful reception to store in a buffer near the actuator. In case of the loss of $\control_{t+\ell}$ at $t+\ell$, $\eta_{t+\ell}$ is applied to the plant if it is already present in buffer, otherwise null control is applied.\label{a:repetitive}	
\end{enumerate}
\begin{figure}
	\begin{adjustbox}{width = \columnwidth}
		\begin{tikzpicture}
		\tikzstyle{block} = [draw, fill=blue!10, rectangle, 
		minimum height=2em, minimum width=1.1cm]
		\tikzstyle{blockyellow} = [draw, fill=yellow!10, rectangle, 
		minimum height=2em, minimum width=1cm]
		\tikzstyle{blockgreen} = [draw, fill=green!10, rectangle, 
		minimum height=2em, minimum width=1.1cm] 
		\tikzstyle{blockred} = [draw, fill=red!10, rectangle, 
		minimum height=2em, minimum width=2cm]  
		\begin{scope}
		\node[coordinate] (0) at (0,0) {};
		%\draw[help lines] (0,0) grid (7,3);
		\node [blockgreen, right=-0.1cm of 0] (controller) {$\offset_t,\gain_t$};
		\node[anchor=north, blue] at (controller.south) {Controller};
		\draw[very thick ,-] (controller.east) -- (2,0);
		\node[coordinate] (1) at (0,1) {};
		\draw[very thick ,->] (1.7,2.2) -- (2.5,2.2);
		\node [blockyellow, right=1.6cm of 1] (packet) {$\control_{t+\ell}$};
		\node [block, above=0cm of packet] (packet) {$(\offset_t)_{\ell m + 1 : m N_r}$};
		\draw[decorate, decoration=switch] (2cm,0cm) -- ++(1cm,0cm);	
		\node [blockred, right=3.5cm of 0] (buffer) {empty};
		\draw[very thick ,-] (3,0) -- (buffer.west);
		\node[anchor=north, blue] at (buffer.south) {Buffer};
		\node [blockred, right=0.5cm of buffer] (actuator) {Actuator};
		\draw[very thick ,-] (buffer.east) -- (actuator.west);
		\draw[->,>=stealth',semithick, blue] (2.6,0.5) arc[radius=0.65, start angle=20, end angle=-45];
		\node[anchor=east, blue] at (3,-0.6) {$\cnoise_{t+\ell}$};
		\end{scope}
		\begin{scope}[shift={(0,-3)}]
		\node[coordinate] (0) at (0,0) {};
		\node [blockgreen, right=-0.1cm of 0] (controller) {$\offset_t,\gain_t$};
		\node[anchor=north, blue] at (controller.south) {Controller};
		\draw[very thick ,-] (controller.east) -- (2,0);
		\node[coordinate] (1) at (0,1) {};
		\draw[very thick ,->] (1.7,1.6) -- (2.5,1.6);
		\node [blockyellow, right=1.6cm of 1] (packet) {$\control_{t+\ell}$};
		\draw[decorate, decoration=switch] (2cm,0cm) -- ++(1cm,0cm);	
		\node [blockred, right=3.5cm of 0] (buffer) {non-empty};
		\draw[very thick ,-] (3,0) -- (buffer.west);
		\node[anchor=north, blue] at (buffer.south) {Buffer};
		\node [blockred, right=0.5cm of buffer] (actuator) {Actuator};
		\draw[very thick ,-] (buffer.east) -- (actuator.west);
		\draw[->,>=stealth',semithick, blue] (2.6,0.5) arc[radius=0.65, start angle=20, end angle=-45];
		\node[anchor=east, blue] at (3,-0.6) {$\cnoise_{t+\ell}$};
		\end{scope}
		\end{tikzpicture}
	\end{adjustbox}	
	\caption{Control channel and buffer at time $t+\ell$ for \ref{a:repetitive}: The blue blocks are transmitted only if the buffer is empty and the yellow blocks are transmitted at each time.}
	\label{fig:policyrep}
\end{figure}
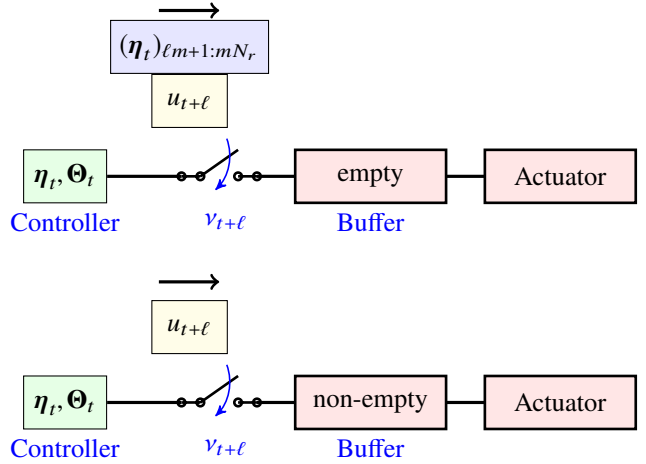	
Since the optimization problem is solved after each recalculation interval, only $N_r \leq N$ blocks of controls in \eqref{e:policy out channel} are applied to the plant, the rest of them are discarded. Therefore, the above protocol transmits $(\offset_t)_{\ell m + 1 : m N_r} = \bmat{ \eta_{t+\ell + 1}\transp & \ldots & \eta_{t+N_r-1}\transp }\transp$ repetitively until the first successful transmission and store them in a buffer at the actuator node. In order to avoid any addition operation at the actuator, $\control_{t+\ell}$ is transmitted at each $t+\ell$ and applied to the plant if successfully received at the actuator otherwise the corresponding $\eta_{t+\ell}$ from the buffer is applied. In a worst case, if the buffer is empty and packets are also lost, then null control is applied to the plant.  		  
The plant input sequence using \ref{a:repetitive} can therefore be represented in compact form as:   
\begin{equation}\label{e:policyrepetitive out}
\control^a_{t:N} \Let \calG \offset_t + \calS \gain_t  \ee(\noisyInnovation_{t:N}), \\
\end{equation}
where \[ \calS \Let \bmat{I_{m}\otimes\cnoise_{t} & & \\ & \ddots & \\ & & I_{m}\otimes \cnoise_{t + \reachindex -1}\\ & & &  I_{m(N-\reachindex)} }, \] and the matrix $\calG$ has $(N \times N)$ blocks in total, each of dimension $m \times m$. For $i = 1,\cdots,N$ and $j = 1,\cdots,N$, the matrix $\calG$ can be given in terms of the blocks $\calG^{(i,j)}$ each of dimension $m \times m$ as follows: 		 
\begin{align}
\calG^{(i,j)} \Let
\begin{cases}
g_{t+i-1} I_m & \text{ if } i=j \leq N_r ,\\ 
I_{m} & \text{ if } i =j > N_r ,\\
0_m \quad  & \text{ otherwise,}			
\end{cases}		
\end{align}
where $g_t = \cnoise_t$, $g_{t+\ell} = g_{t+\ell-1} + (1- g_{t+\ell-1})\cnoise_{t+\ell}  $, and $\gain_t$ and $\ee(\noisyInnovation_{t:N})$ are as defined in \eqref{e:policy out channel}. The term $g_{t+\ell}$ captures the effect of \ref{a:repetitive}. 
Note that \ref{a:repetitive} requires storage at the actuator, but no advanced computation capacity.

\section{Tractability}\label{s:tractability}
In this section we present a tractable surrogate of the optimal control problem \eqref{e:opt control problem} under the system setup discussed in \secref{s:setup}. 
Let us define $\mu_{\calG} \Let \EE[\calG]$, $\Sigma_{\calG} \Let \EE \left[ \calG\transp \alpha \calG \right]$, $\alpha \Let \calB \transp \calQ \calB + \calR$, $\mu_{\calS} \Let \EE[\calS]$, $\Sigma_{\calS} \Let \EE[\calS\transp \alpha \calS]$, $\Sigma_{\calG\calS} \Let \EE[\calG\transp \alpha \calS ]$, $\Pi_{\meas_t} \Let \ee_0(\noisyInnovation_t) \ee_0(\noisyInnovation_t) \transp $, $\Sigma_{\ee} \Let \EE\Bigl[\ee^{\prime}(\noisyInnovation_{t+1:N-1})\ee^{\prime}(\noisyInnovation_{t+1:N-1})\transp \Bigr]$, \\
$\Sigma_{\ee^{\prime}\wnoise} \Let \EE \Biggl[\ee^{\prime}(\noisyInnovation_{t+1:N-1})\wnoise_{t:N} \transp \Biggr]$, $\Sigma_{e\ee^{\prime}} \Let \EE \Bigl[ \ee^{\prime} (\noisyInnovation_{t+1:N-1})e_t \transp \Bigr]$, $\gain_t^{(:,t)} \Let \bmat{\theta_{0,t}\transp & \theta_{1,t} \transp & \cdots & \theta_{N-1,t}\transp }\transp $, $\calQ_{\calA} \Let \calA \transp \calQ \calB$ and $\calQ_{\calD} \Let \calD \transp \calQ \calB$.  
We have the following Lemma:		
\begin{lem}\label{lem:objective}
	The objective function \eqref{e:cost_compact} can be represented in terms of the decision variables as follows:
	\begin{equation} \label{e:obj out channel}
	\begin{aligned}
	& V_t^{\prime} = \offset_t\transp \Sigma_{\calG} \offset_t + \trace(\Sigma_{\calS} \gain_t^{(:,t)} \Pi_{\meas_t} (\gain_t^{(:,t)})\transp) + \trace (\Sigma_{\calS} \gain_t^{\prime} \Sigma_{\ee} (\gain_t^{\prime})\transp) \\
	& +  2(\offset_t \transp \Sigma_{\calG\calS} + \stest_t \transp \calQ_{\calA} \mu_{\calS} ) \gain_t^{(:,t)}\ee_0(\noisyInnovation_t) +2\stest_t \transp \calQ_{\calA}\mu_{\calG}\offset_t \\ 
	&+ 2\trace (\calQ_{\calD} \mu_{\calS} \gain_t^{\prime} \Sigma_{\ee^{\prime}\wnoise}) + 2\trace \Bigl( \calQ_{\calA} \mu_{\calS} \gain_t^{\prime} \Sigma_{e \ee^{\prime}} \Bigr) .
	\end{aligned}
	\end{equation}
\end{lem}
A proof of the Lemma \ref{lem:objective} is given in the appendix. At each optimization time $t$, the above objective function is updated by substituting new information $\noisyInnovation_{t}$ and $\stest_t$. The matrices $\mu_{\calG}$, $\Sigma_{\calG}, \mu_{\calS}, \Sigma_{\calS}, \Sigma_{\calG\calS}, \Sigma_{\ee}, \Sigma_{\ee^{\prime}\wnoise}, \Sigma_{e\ee^{\prime}} $ incorporate all elements of the problem setup (\secref{s:setup}). These matrices are computed offline by using Monte-Carlo simulations to reduce the burden of online computation. We have the following result:	

\begin{prop}\label{th:main}
	For every time $t= 0, N_r, 2N_r, \cdots$, the optimal control problem \eqref{e:opt control problem} can be written as the following convex quadratic, (globally) feasible program:
	\begin{align}\label{e:main program}
	\minimize_{\offset_t,\gain_t} & \quad  \eqref{e:obj out channel} \notag \\  %
	\sbjto & \quad \eqref{e:decisionboundsingle}
	\end{align}
\end{prop}
\begin{pf}
	The objective function \eqref{e:obj out channel} is convex quadratic in decision variables $\offset_t$ and $\gain_t$, and the constraint \eqref{e:decisionboundsingle} is a convex affine function of the decision variables. Since \eqref{e:decisionboundsingle} does not depend on $\st_t$, the optimization program \eqref{e:main program} is feasible for all $\st_t \in \R^d$ for all $t \in \Nz$.
\end{pf}	
\section{Stability}\label{s:stability}
In this section, we show that the system setup discussed in \secref{s:setup} leads to the mean square boundedness of the controlled states if carefully designed stability constraints are also included in the optimization program \eqref{e:main program}. 
Let us first represent the estimator process recursion \eqref{e:est} in terms of the matrix pair $(A,B)$ as follows:
\begin{equation}\label{e:estEquation}
\stest_{t+1} = \A \stest_t + \B \control_t^a + \tilde{\wnoise}_t,
\end{equation}
where $\tilde{\wnoise}_t \Let \snoise_{t+1} (\A \estError_{t} + \hat{\wnoise}_{t})$. 
Note that a Lyapunov stable system matrix $\A$ can be decomposed into a Schur stable component $\A_s$ and an orthogonal component $\A_o$ as:
\begin{equation}\label{e:decomposed st}
\bmat{{\stest}_{t+1}^o \\ {\stest}_{t+1}^s} = \bmat{\Aortho {\stest}_{t}^o \\ \Aschur {\stest}_{t}^s} + \bmat{\Bortho  \\ \Bschur} \control_t^a + \bmat{\tilde{\wnoise}_t^o \\ \tilde{\wnoise}_t^s},
\end{equation}
where ${\stest}_t^s \in \R^{\dschur}$ , ${\stest}_t^o \in \R^{\dortho}$, and $d = \dortho + \dschur $. Let us define the reachability matrix
\begin{equation}\label{e:reachabilityMatrix} 
\reachab_{\reachindex}(\A,\B) \Let \bmat{\A^{\reachindex  -1}\B & \ldots & \A \B & \B}. 
\end{equation}
By the stabilizability of $(\A,\B)$, there exists an integer $\reachindex$ such that the reachability matrix $ \reachab_{\reachindex}(\Aortho,\Bortho) $ has full row rank. The integer $\reachindex$ is called reachability index of the matrix pair $(\Aortho,\Bortho)$. The reachability index is important in our approach. We can choose $N_r \geq \reachindex$ but for simplicity, in our subsequent analysis we use $N_r = \reachindex$.
 We consider the orthogonal component of the $\reachindex$-subsampled process of \eqref{e:estEquation}, which is given by   
	\begin{equation} \label{e:kappaEstimator}
	\stest^o_{\reachindex(t+1)} = \Aortho^{\reachindex}\stest^o_{\reachindex t} + \reachab_{\reachindex}(\Aortho,\Bortho)\control_{\reachindex t:\reachindex}^a + \reachab_{\reachindex}(\Aortho,I)\tilde{\wnoise}^o_{\reachindex t: \reachindex}.
	\end{equation}
Let us define $z_t \Let \left( (\Aortho^{\reachindex t})\transp \stest_{\reachindex t}^o \right)$ then the process $(z_t)_{t \in \Nz}$ can be considered as a d-dimensional random walk with recursion
\begin{equation}\label{e:z_difference}
z_{t+1} = z_t + (\Aortho^{\reachindex(t+1)})\transp\left( \reachab_{\reachindex}(\Aortho, \Bortho)\control_{\reachindex t: \reachindex}^a + \reachab_{\reachindex}(\Aortho, I)\tilde{\wnoise}_{\reachindex t: \reachindex}^o\right).
\end{equation}
We present the following lemma:  
\begin{lem} \label{lem:ortho stable intermittent}
	Consider the recursions \eqref{e:estEquation} and \eqref{e:z_difference}.  Suppose that $\control_t^a$ is constrained in the set $\controlset$ for each $t$ and that there exist $a, r > 0$ such that for $ j=1,2,\cdots ,d_o$, the following conditions hold
%	\eqref{e:drift general} hold,
	\begin{subequations}\label{e:drift general}
	\begin{align}
&	\EE_{\XX_{\reachindex t}} \left[ \left(z_{t+1}\right)^{(j)} - \left(z_t\right)^{(j)} \right]  \leq -a
	\text{ whenever }  \left( z_t \right)^{(j)} > r, \label{e:drift general 1}\\
&	\EE_{\XX_{\reachindex t}} \left[ \left(z_{t+1}\right)^{(j)} - \left(z_t\right)^{(j)} \right]  \geq a  
	\text{ whenever }  \left( z_t  \right)^{(j)} < -r,  \label{e:drift general 2} \\
&	\EE \left[ \abs{\left(z_{t+1}\right)^{(j)} - \left(z_t\right)^{(j)} }^4 \mid z_0^{(j)}, \ldots, z_t^{(j)} \right]  \leq M  
		\text{ for all }  t.  \label{e:drift general 3}  
	\end{align}
\end{subequations}
	Then there exists $\bar{\gamma}> 0$ such that $\EE_{\XX_0} \left[ \norm{\stest_t}^2 \right] \leq \bar{\gamma}$ for all $t \geq 0$. Moreover, there exists a $\reachindex$-history dependent class of policies $\control_{\reachindex t : \reachindex}$ of the form \eqref{e:policy out channel} such that \eqref{e:drift general} and \eqref{e:controlset} are satisfied for $\zeta \in \left]0, \frac{\authority}{\sqrt{d_o}\sigma_{1}\left(\reachab_{\reachindex}(\Aortho,\Bortho)^{\dagger}\right)} \right] $ and $a = \zeta p_c$ under the transmission protocol \ref{a:repetitive}.
	Furthermore, for $t=0,\reachindex, 2\reachindex, \ldots$, the conditions \eqref{e:drift general} are equivalent to the following conditions:
	\begin{subequations} \label{e:drift intermittent}
		\begin{align}
	&	\Bigl( (\Aortho^{t+ \reachindex })\transp\reachab_{\reachindex}(\Aortho, \Bortho)\EE_{\XX_t} [\control_{t:\reachindex}] \Bigr)^{(j)} \leq -\zeta \notag  \\
		& \quad \text{ whenever }  \left( (\Aortho^t) \transp \stest_{t}^o \right)^{(j)} > r, \label{e:drift1 intermittent} \\	
	&	\Bigl( (\Aortho^{t+  \reachindex})\transp\reachab_{\reachindex}(\Aortho, \Bortho)\EE_{\XX_t} [\control_{t:\reachindex}]\Bigr)^{(j)} \geq \zeta \notag \\
	& 	\quad  \text{ whenever }  \left( (\Aortho^t) \transp \stest_{t}^o \right)^{(j)} < -r. \label{e:drift2 intermittent}
		\end{align}
	\end{subequations}
\end{lem}
A proof of Lemma \ref{lem:ortho stable intermittent} is given in the appendix.
To embed the drift conditions \eqref{e:drift intermittent} in a tractable optimization program, we consider the first $\reachindex$ blocks of \eqref{e:policy out channel} 
\begin{equation} \label{e:kappa blocks control}
\control_{t:\reachindex} = (\offset_t)_{1:\reachindex m} + (\gain_t)_{1:\reachindex m } \ee(\noisyInnovation_{t:N})
\end{equation} 
for $t=0, \reachindex, \cdots$, and substitute them in \eqref{e:drift intermittent}. We get the following stability constraints:
\begin{subequations}\label{e:stability constraint}
	\begin{align}
	& \quad \Biggl( (\Aortho^{t+\reachindex})\transp\reachab_{\reachindex}(\Aortho, \Bortho) \left( (\offset_t)_{1: \reachindex m} + (\gain_t^{(:,t)})_{1:\reachindex m}\ee_0(\noisyInnovation_t) \right) \Biggr)^{(j)} \leq -\zeta \nonumber  \\
	& \quad \quad \text{ whenever } \left( (\Aortho^{ t})\transp\stest^o_{t}  \right)^{(j)} > r, \label{e:stability constraint 1}\\
	& \quad \Biggl( (\Aortho^{t+\reachindex})\transp\reachab_{\reachindex}(\Aortho, \Bortho)\left( (\offset_t)_{1: \reachindex m} + (\gain_t^{(:,t)})_{1:\reachindex m}\ee_0(\noisyInnovation_t) \right) \Biggr)^{(j)} \geq \zeta \nonumber \\ 
	& \quad \quad \text{ whenever } \left( (\Aortho^{ t})\transp \stest^o_{t} \right)^{(j)} < -r , \label{e:stability constraint 2} 
	\end{align}
\end{subequations}
where $r, \zeta$ and $j$ are as in \eqref{e:drift intermittent}. We have the following result:

\begin{thm}\label{th:stbl}
	Let the stability constraints \eqref{e:stability constraint} be included in the optimization program of Theorem \ref{th:main} and controls be generated by successively solving the underlying optimization program. Then the application of these controls ensures mean-square boundedness of the states of \eqref{e:system} for any bound $\authority > 0$ and transmission probabilities $p_s,  p_c > 0$ . 
\end{thm}
A proof of Theorem \ref{th:stbl} is given in the appendix.

\section{Numerical Experiments}\label{s:experiment}

%====================================
\begin{figure}[t]
	\centering
	\begin{adjustbox}{width=\columnwidth}
		\begin{tikzpicture}		
		\begin{axis}[%
		width=3.028in,
		height=2.754in,
		at={(2.509in,1.111in)},
		scale only axis,
		xmin=1,
		xmax=5,
		xtick={1,2,3,4,5},
		xticklabels={{2},{3},{4},{5},{10}},
		xlabel={$\authority$},
		xmajorgrids,
		ymin=600,
		ymax=1200,
		ylabel={Empirical MSB},
		ymajorgrids,
		axis background/.style={fill=white},
		legend style={at={(0.423,0.562)},anchor=south west,legend cell align=left,align=left,draw=white!15!black}
		]
		\addplot [color=red,solid,line width=2pt]
		table[row sep=crcr]{%
			1	1170.96160052718\\
			2	856.799512445458\\
			3	733.914388537433\\
			4	673.594949506186\\
			5	628.592066078147\\
		};
		\addlegendentry{$\gain_t$ as in \eqref{e:gain}};
		
		\addplot [color=green,dashed,line width=2pt]
		table[row sep=crcr]{%
			1	1196.12242203644\\
			2	886.72043804782\\
			3	771.326230670421\\
			4	714.832485424553\\
			5	685.78954428522\\
		};
		\addlegendentry{$\gain_t = \zeros$};
		\addplot [color=blue,dotted,line width=2pt]
		table[row sep=crcr]{%
			1	1172.38882443538\\
			2	861.885113633637\\
			3	738.39198463386\\
			4	681.078550433424\\
			5	645.587090856976\\
		};
		\addlegendentry{$\gain_t = \tilde{\gain}_t$ as in \eqref{e:gain2}};	
		\addplot [color=black, mark = o]
		table[row sep=crcr]{%
			1	1093.0370266828\\
			2	822.509585662715\\
			3	717.751433720204\\
			4	670.875619046365\\
			5	626.466316775453\\
		};
		\addlegendentry{without \eqref{e:stability constraint}}	
		\end{axis}
		
		\end{tikzpicture}%
	\end{adjustbox}
	\caption{Empirical MSB when $\authority$ varies in the set $\{2,3,4,5,10 \}$ while $p_c=p_s = 0.8$ remain fixed. }
	\label{Fig:varyingAuthority}
\end{figure}
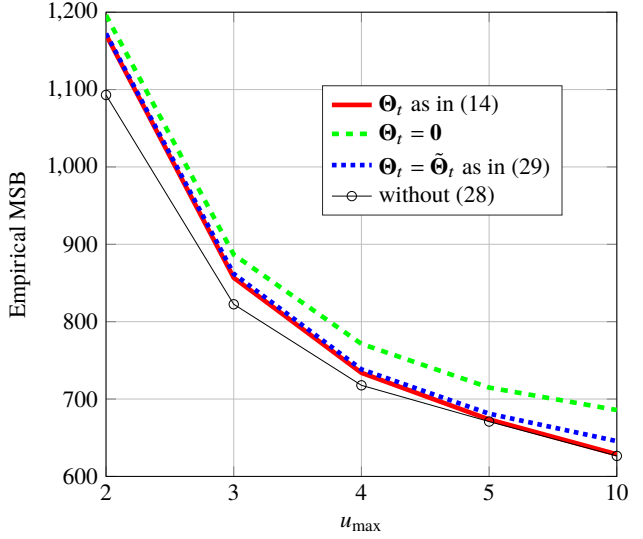

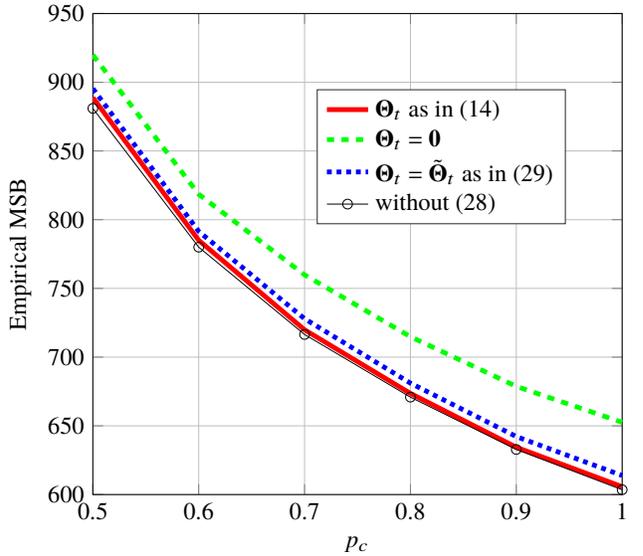
\begin{figure}[t]
	\centering
	\begin{adjustbox}{width=\columnwidth}
		\begin{tikzpicture}		
		\begin{axis}[%
		width=3.028in,
		height=2.754in,
		at={(2.509in,1.111in)},
		scale only axis,
		xmin=1,
		xmax=6,
		xtick={1,2,3,4,5,6},
		xticklabels={{0.5},{0.6},{0.7},{0.8},{0.9},{1}},
		xlabel={$p_c$},
		xmajorgrids,
		ymin=600,
		ymax=950,
		ylabel={Empirical MSB},
		ymajorgrids,
		axis background/.style={fill=white},
		legend style={at={(0.423,0.562)},anchor=south west,legend cell align=left,align=left,draw=white!15!black}
		]
		\addplot [color=red,solid,line width=2pt]
		table[row sep=crcr]{%
		1	888.821141200948\\
		2	785.005298740941\\
		3	719.878334351184\\
		4	673.594949506186\\
		5	634.418354233448\\
		6	605.69204022682\\
		};
		\addlegendentry{$\gain_t$ as in \eqref{e:gain}};
		
		\addplot [color=green,dashed,line width=2pt]
		table[row sep=crcr]{%
			1	919.833678049158\\
			2	818.510031015265\\
			3	759.873349945105\\
			4	714.832485424553\\
			5	678.72211654104\\
			6	652.565242124556\\
		};
		\addlegendentry{$\gain_t = \zeros$};
		\addplot [color=blue,dotted,line width=2pt]
		table[row sep=crcr]{%
			1	895.413997636264\\
			2	791.495083270434\\
			3	727.928882033558\\
			4	681.078550433424\\
			5	642.208811807097\\
			6	613.861500485622\\
		};
		\addlegendentry{$\gain_t= \tilde{\gain}_t$ as in \eqref{e:gain2}};	
		\addplot [color=black, mark=o]
		table[row sep=crcr]{%
			1	880.96609723598\\
			2	779.991761647779\\
			3	716.484268996089\\
			4	670.875619046365\\
			5	632.819068509409\\
			6	603.502045804726\\
		};
		\addlegendentry{without \eqref{e:stability constraint}}	
		\end{axis}
		\end{tikzpicture}%
	\end{adjustbox}
	\caption{Empirical MSB when $p_c$ varies from $0.5$ to $1$ while $\authority = 5$ and $p_s = 0.8$ remain fixed. }
	\label{Fig:varyingCnoise}
\end{figure}

\begin{figure}[t]
	\centering
	\begin{adjustbox}{width=\columnwidth}
		\begin{tikzpicture}		
		\begin{axis}[%
		width=3.028in,
		height=2.754in,
		at={(2.509in,1.111in)},
		scale only axis,
		xmin=1,
		xmax=6,
		xtick={1,2,3,4,5,6},
		xticklabels={{0.5},{0.6},{0.7},{0.8},{0.9},{1}},
		xlabel={$p_s$},
		xmajorgrids,
		ymin=620,
		ymax=820,
		ylabel={Empirical MSB},
		ymajorgrids,
		axis background/.style={fill=white},
		legend style={at={(0.423,0.562)},anchor=south west,legend cell align=left,align=left,draw=white!15!black}
		]
		\addplot [color=red,solid,line width=2pt]
		table[row sep=crcr]{%
			1	793.338898311148\\
			2	730.706536704229\\
			3	693.795385432696\\
			4	673.594951282963\\
			5	653.485938362923\\
			6	634.649865879644\\
		};
		\addlegendentry{$\gain_t$ as in \eqref{e:gain}};
		
		\addplot [color=green,dashed,line width=2pt]
		table[row sep=crcr]{%
			1	819.068305028279\\
			2	760.886208211163\\
			3	729.202244964196\\
			4	714.832485450093\\
			5	702.14305798105\\
			6	687.989779721633\\
		};
		\addlegendentry{$\gain_t = \zeros$};
		\addplot [color=blue,dotted,line width=2pt]
		table[row sep=crcr]{%
			1	798.062108896244\\
			2	737.511012435875\\
			3	702.143419446884\\
			4	681.078550433424\\
			5	659.728580120119\\
			6	644.080977214542\\
		};
		\addlegendentry{$\gain_t = \tilde{\gain}_t$ as in \eqref{e:gain2}};
		\addplot [color=black, mark = o]
		table[row sep=crcr]{%
			1	789.394083458429\\
			2	727.511080467208\\
			3	691.970897812241\\
			4	670.875619046365\\
			5	650.163060713868\\
			6	632.996409612861\\
		};	
	    \addlegendentry{without \eqref{e:stability constraint}}	
		\end{axis}
		\end{tikzpicture}%
	\end{adjustbox}
	\caption{Empirical MSB when $p_s$ varies from $0.5$ to $1$ while $\authority = 5$ and $p_c = 0.8$ remain fixed. }
	\label{Fig:varyingSnoise}
\end{figure}
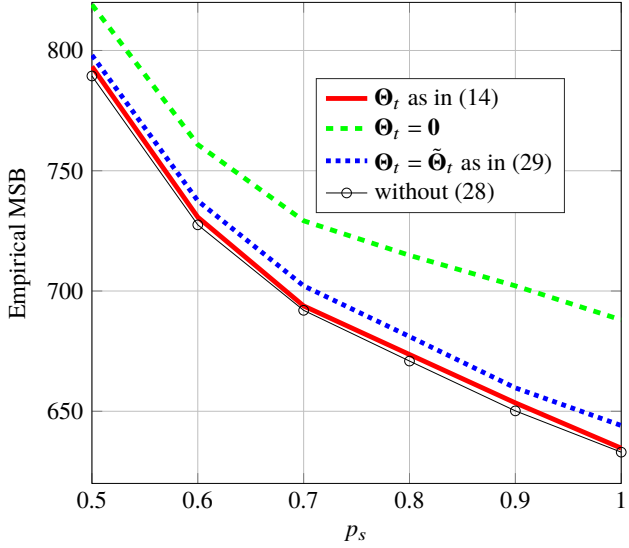

%======================================================================
In this section, we present numerical experiments and record the empirical mean of the quantities of interest to illustrate our results by taking averages over 1000 sample paths for 120 time steps. We consider the four dimensional stochastic LTI system \eqref{e:system} with matrices taken from \cite{ref:Hokayem-12},
which can be written in the form of \eqref{e:decomposed st} with $\Aschur = 0.9$, $\Bschur = 0 $, \[\Aortho = \bmat{1 & 0 & 0\\0 & 0 & -1\\ 0 & 1 & 0} , \text{ and } \Bortho = \bmat{1 \\ 0 \\ 1}. \]
The reachability matrix (defined in \eqref{e:reachabilityMatrix})
\[ \reachab_{3}(\Aortho, \Bortho) = \bmat{
1 &	1&	1 \\
0&	-1&	0 \\
-1&	0	&1} 
\]
has full row rank for $\reachindex  = 3$.
We repeatedly solved a finite-horizon CSOCP reported in Theorem \ref{th:stbl} corresponding to states and control weights \( Q = I_4 , Q_f = I_4, R =1 \), the optimization horizon, \(N=5\), recalculation interval \(N_r = \reachindex = 3 \) and simulation data $\stinit \sim N(0,I_4)$, $\wnoise_t \sim N(0,10I_4)$, $\mnoise_t \sim N(0,10I_4)$. We selected the nonlinear bounded term $\ee(\cdot)$ in our policy, similar to \cite{PDQ-LCSS}, to be a vector of scalar sigmoidal functions $ \varphi(\xi)=\frac{1- e^{-\xi}}{1+ e^{-\xi}} $ applied to each coordinate of the received innovation sequence. We use a MATLAB-based software package YALMIP \cite{ref:lofberg-04} and a solver SDPT3-4.0 \cite{ref:toh-06} to solve the underlying optimization programs.  
\par
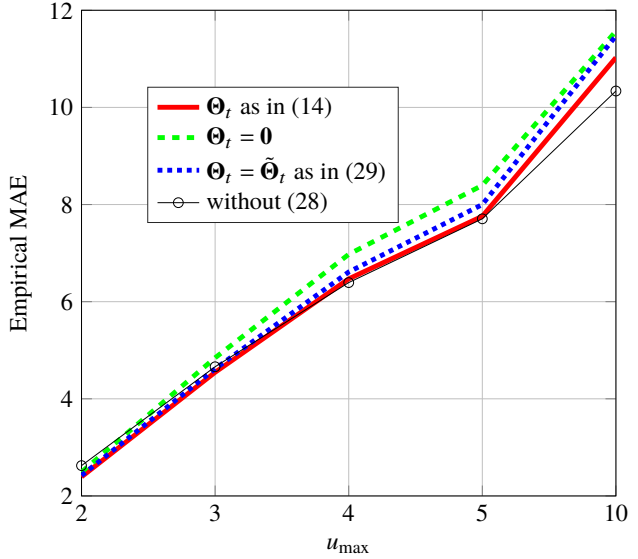
\begin{figure}[t]
	\centering
	\begin{adjustbox}{width=\columnwidth}
		\begin{tikzpicture}		
		\begin{axis}[%
		width=3.028in,
		height=2.754in,
		at={(2.509in,1.111in)},
		scale only axis,
		xmin=1,
		xmax=5,
		xtick={1,2,3,4,5},
		xticklabels={{2},{3},{4},{5},{10}},
		xlabel={$\authority$},
		xmajorgrids,
		ymin=2,
		ymax=12,
		ylabel={Empirical MAE},
		ymajorgrids,
		axis background/.style={fill=white},
		legend style={at={(0.123,0.562)},anchor=south west,legend cell align=left,align=left,draw=white!15!black}
		]
		\addplot [color=red,solid,line width=2pt]
		table[row sep=crcr]{%
			1	2.39355202050923\\
			2	4.54588680068135\\
			3	6.46980237894085\\
			4	7.75860955948248\\
			5	11.0250514763607\\
		};
		\addlegendentry{$\gain_t$ as in \eqref{e:gain}};
		
		\addplot [color=green,dashed,line width=2pt]
		table[row sep=crcr]{%
		1	2.49083502178213\\
		2	4.8445629868086\\
		3	6.97485763820439\\
		4	8.3954607165089\\
		5	11.5694121489621\\
		};
		\addlegendentry{$\gain_t = \zeros$};
		\addplot [color=blue,dotted,line width=2pt]
		table[row sep=crcr]{%
		1	2.42129411625368\\
		2	4.61852922111149\\
		3	6.61164514181819\\
		4	7.99986698930331\\
		5	11.4728169655117\\
		};
		\addlegendentry{$\gain_t = \tilde{\gain}_t$ as in \eqref{e:gain2}};
		\addplot [color=black, mark = o]
		table[row sep=crcr]{%
			1	2.62471625539073\\
			2	4.65985871961617\\
			3	6.39309902519719\\
			4	7.70629743045316\\
			5	10.3376967631189\\
		};
		\addlegendentry{without \eqref{e:stability constraint}}		
		\end{axis}
		\end{tikzpicture}%
	\end{adjustbox}
	\caption{Empirical MAE per stage when $\authority$ varies in the set $\{2,3,4,5,10 \}$ while $p_c=p_s = 0.8$ remain fixed. }
	\label{Fig:varyingUmaxEnergy}
\end{figure}

 We compare the present approach with a simplified version of our main result by setting $\gain_t = \zeros$ in \eqref{e:policy out channel} and using a modified version of $\gain_t$ as giveen below:
\begin{equation} \label{e:gain2}
\tilde{\gain}_t = \bmat{ \theta_{0, t} & \zeros & \cdots & \zeros & \zeros \\ \zeros & \theta_{1, t+1}  & \cdots & \zeros & \zeros \\ \vdots & \vdots & \vdots & \vdots & \vdots\\ \zeros & \zeros & \cdots & \zeros & \theta_{N-1, t+N-1} }.
\end{equation}
  In particular, $\gain_t = \zeros$ represents the optimization over an open-loop control sequence and $\gain_t = \tilde{\gain}_t$ represents the optimization over only one causal feedback term. 
	In the above three cases, we have used stability constraints \eqref{e:stability constraint} as mentioned in Theorem \ref{th:stbl}. In the fourth case, we have removed stability constraints \eqref{e:stability constraint} and simulated with the same simulation data.
Our observations from numerical experiments are listed below. 
\begin{enumerate}[label={\rm (O\arabic*)},leftmargin = *, nosep, start = 1]
	\item In Fig.\ \ref{Fig:varyingAuthority}, we plot the empirical mean square bound (MSB) with respect to  $\authority$ picked from the set $\{2,3,4,5,10\}$ while $p_s = p_c = 0.8$ remain fixed. The empirical MSB decreases with the increase in $\authority$. The empirical MSB is less when $\gain_t$ is chosen according to \eqref{e:gain} and it further decreases when we remove the constraints \eqref{e:stability constraint}. The difference due to the removal of \eqref{e:stability constraint} is more for the lower values of $\authority$ and it vanishes when we further increase $\authority$. \label{o:umsb}
	\item In Fig. \ref{Fig:varyingCnoise},  we plot the empirical MSB with respect to  $p_c$ picked from the set $\{i/10 \mid i = 5, \ldots 10 \}$ while $p_s = 0.8$, $\authority = 5$ remain fixed. The empirical MSB decreases with increase in $p_c$ with higher slope at smaller $p_c$. The empirical MSB for $\gain_t = \zeros$ is the highest of all considered cases and the difference increases with the increase in $p_c$. The lowest empirical MSB's are achieved in the absence of \eqref{e:stability constraint}. \label{o:cmsb}
	\item In Fig. \ref{Fig:varyingSnoise}, we plot the empirical MSB with respect to $p_s$ picked from the set $\{i/10 \mid i = 5, \ldots 10 \}$ while $p_c = 0.8$, $\authority = 5$ remain fixed. With increase in $p_s$, we observe decrease in empirical MSB. In this case also, the lowest empirical MSB's are achieved in the absence of \eqref{e:stability constraint}. \label{o:smsb}
	\item In Fig. \ref{Fig:varyingUmaxEnergy},  we plot the empirical mean of actuator energy (MAE) per stage with respect to  $\authority$ by fixing the parameters as in the observation \ref{o:umsb}. The empirical MAE per stage increases with the increase in $\authority$. The empirical MAE per stage in the absence of \eqref{e:stability constraint} is the largest (slightly) for $\authority = 2$  and the lowest for $\authority = 10$ of all considered cases.
	\item In Fig. \ref{Fig:varyingCnoiseEnergy}, we plot the empirical MAE per stage with respect to $p_c$ by fixing the parameters as in the observation \ref{o:cmsb}. The empirical MAE per stage increases with the increase in $p_c$. In this observation, the empirical MAE per stage is the lowest in the absence of \eqref{e:stability constraint}.
	\item In Fig. \ref{Fig:varyingSnoiseEnergy},  we plot the empirical MAE per stage with respect to  $p_s$ by fixing the parameters as in the observation \ref{o:smsb}. The empirical MAE per stage does not vary much with the increase in $p_s$. The lowest empirical MAE per stage is obtained in the absence of \eqref{e:stability constraint}.
\end{enumerate}
The proposed class of policies performs better than open loop control sequence (which is obtained by substituting $\gain_t = \zeros$) in terms of empirical MSB and empirical MAE. Moreover, numerical experiments demonstrate mean square boundedness of controlled states for all considered cases. The optimization over $\tilde{\gain}_t$ in place of $\gain_t$ results in significant reduction in the number of decision variables (please also see \cite[Remark 1]{ref:ChaHokLyg-11}) with minimal increase in MSB and MAE. 
\par
\begin{figure}
	\begin{adjustbox}{width=\columnwidth}
%		\input{Umax_runtime.tex}
	\begin{tikzpicture}
	
	\begin{axis}[%
	width=5.52in,
	height=4.556in,
	at={(0.758in,0.487in)},
	scale only axis,
	bar shift auto,
	xmin=0.511111111111111,
	xmax=5.48888888888889,
	xtick={1,2,3,4,5},
	xticklabels={{2},{3},{4},{5},{10}},
	xlabel style={font=\Large},
	xlabel={$\authority$},
	ymin=-10,
	ymax=85,
	ylabel style={font=\Large},
	ylabel={Percentage difference in solver-time},
	axis background/.style={fill=white},
	xmajorgrids,
	ymajorgrids,
	legend style={legend cell align=left, align=left, draw=white!15!black, font = \Large}
	]
	\addplot[ybar, bar width=0.178, fill=green, draw=black, area legend] table[row sep=crcr] {%
		1	64.3583649354721\\
		2	63.1451952404974\\
		3	63.0021220732892\\
		4	64.2501933284617\\
		5	64.5265095930141\\
	};
	\addplot[forget plot, color=white!15!black] table[row sep=crcr] {%
		0.511111111111111	0\\
		5.48888888888889	0\\
	};
	\addlegendentry{$\gain_t = \zeros$}
	
	\addplot[ybar, bar width=0.178, fill=blue, draw=black, area legend] table[row sep=crcr] {%
		1	43.6346285399253\\
		2	41.8317894693203\\
		3	41.880438329666\\
		4	41.6348651790351\\
		5	43.0938770264141\\
	};
	\addplot[forget plot, color=white!15!black] table[row sep=crcr] {%
		0.511111111111111	0\\
		5.48888888888889	0\\
	};
	\addlegendentry{$\gain_t = \tilde{\gain}_t$ as in \eqref{e:gain2}}
	
	\addplot[ybar, bar width=0.178, fill=black, draw=black, area legend] table[row sep=crcr] {%
		1	3.24348012609616\\
		2	0.34183318153698\\
		3	1.09176434867417\\
		4	1.28508386952537\\
		5	-0.39298548516609\\
	};
	\addplot[forget plot, color=white!15!black] table[row sep=crcr] {%
		0.511111111111111	0\\
		5.48888888888889	0\\
	};
	\addlegendentry{without \eqref{e:stability constraint}}
	
	\end{axis}
	
	\begin{axis}[%
	width=5.832in,
	height=4.371in,
	at={(0in,0in)},
	scale only axis,
	xmin=0,
	xmax=1,
	ymin=0,
	ymax=1,
	axis line style={draw=none},
	ticks=none,
	axis x line*=bottom,
	axis y line*=left,
	legend style={legend cell align=left, align=left, draw=white!15!black}
	]
	\end{axis}
	\end{tikzpicture}%		
	\end{adjustbox}	
	\caption{Percentage difference in solver-time with respect to the Theorem \ref{th:stbl} when $\authority$ varies and $p_c = p_s = 0.8$ remain fixed.}
	\label{fig:runtim_authority}
\end{figure}
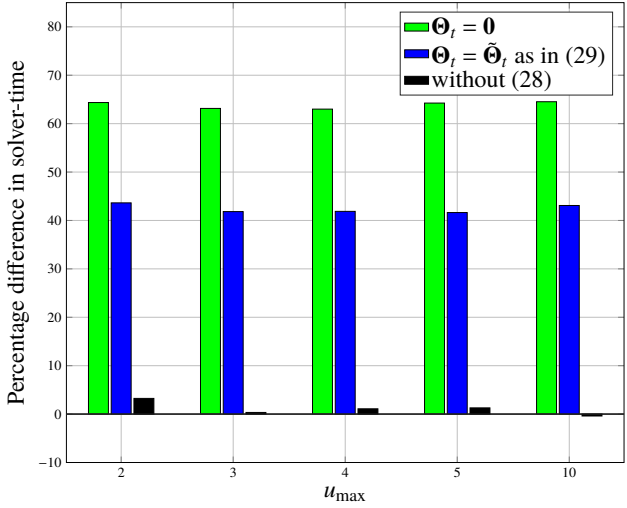

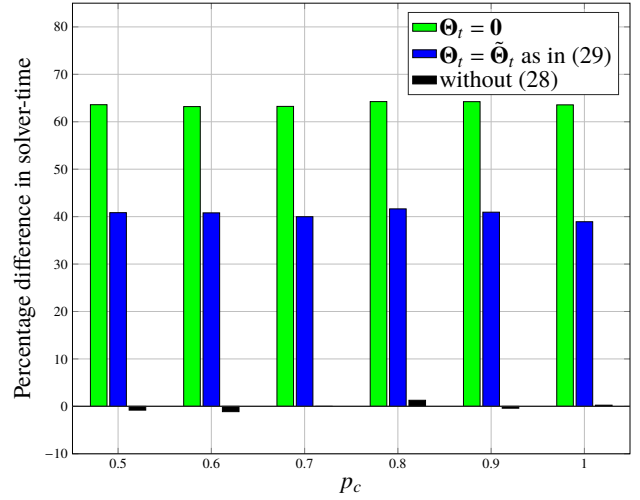
\begin{figure}
	\begin{adjustbox}{width=\columnwidth}
\begin{tikzpicture}

\begin{axis}[%
width=5.574in,
height=4.509in,
at={(0.767in,0.489in)},
scale only axis,
bar shift auto,
xmin=0.511111111111111,
xmax=6.48888888888889,
xtick={1,2,3,4,5,6},
xticklabels={{0.5},{0.6},{0.7},{0.8},{0.9},{1}},
xlabel style={font=\Large},
xlabel={$p_c$},
ymin=-10,
ymax=85,
ylabel style={font=\Large},
ylabel={Percentage difference in solver-time},
axis background/.style={fill=white},
xmajorgrids,
ymajorgrids,
%legend style={anchor=south west,legend cell align=left,align=left,draw=white!15!black}
%legent style = {font = \Large}
legend style={legend cell align=left, align=left, draw=white!15!black, font = \Large}
]
\addplot[ybar, bar width=0.178, fill=green, draw=black, area legend] table[row sep=crcr] {%
	1	63.5981178521311\\
	2	63.197980207126\\
	3	63.2408705691927\\
	4	64.2501933284617\\
	5	64.2369922650602\\
	6	63.5628896104926\\
};
\addplot[forget plot, color=white!15!black] table[row sep=crcr] {%
	0.511111111111111	0\\
	6.48888888888889	0\\
};
\addlegendentry{$\gain_t = \zeros$}

\addplot[ybar, bar width=0.178, fill=blue, draw=black, area legend] table[row sep=crcr] {%
	1	40.8362129618458\\
	2	40.7944568673361\\
	3	39.9806751454213\\
	4	41.6348651790351\\
	5	40.937982023429\\
	6	38.9249935398901\\
};
\addplot[forget plot, color=white!15!black] table[row sep=crcr] {%
	0.511111111111111	0\\
	6.48888888888889	0\\
};
\addlegendentry{$\gain_t = \tilde{\gain}_t$ as in \eqref{e:gain2}}

\addplot[ybar, bar width=0.178, fill=black, draw=black, area legend] table[row sep=crcr] {%
	1	-0.794211404051654\\
	2	-1.1169527625708\\
	3	-0.0186475750254239\\
	4	1.28508386952537\\
	5	-0.385775195542083\\
	6	0.250915361037431\\
};
\addplot[forget plot, color=white!15!black] table[row sep=crcr] {%
	0.511111111111111	0\\
	6.48888888888889	0\\
};
\addlegendentry{without \eqref{e:stability constraint}}

\end{axis}

\begin{axis}[%
width=5.902in,
height=4.322in,
at={(0in,0in)},
scale only axis,
xmin=0,
xmax=1,
ymin=0,
ymax=1,
axis line style={draw=none},
ticks=none,
axis x line*=bottom,
axis y line*=left,
legend style={legend cell align=left, align=left, draw=white!15!black}
]
\end{axis}
\end{tikzpicture}%	
	\end{adjustbox}	
	\caption{Percentage difference in solver-time with respect to the Theorem \ref{th:stbl} when $p_c$ varies and $p_s = 0.8$, $\authority = 5$ remain fixed.}
	\label{fig:runtime_cnoise}
\end{figure}

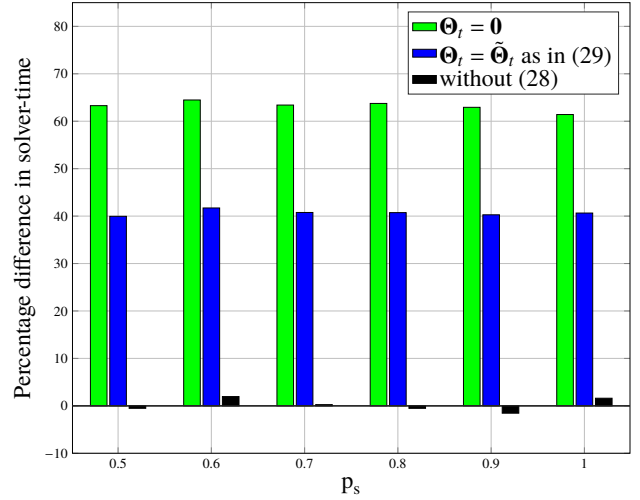
\begin{figure}
	\begin{adjustbox}{width=\columnwidth}
%%		\input{runtime_snoise.tex}	
\begin{tikzpicture}

\begin{axis}[%
width=5.574in,
height=4.509in,
at={(0.767in,0.489in)},
scale only axis,
bar shift auto,
xmin=0.511111111111111,
xmax=6.48888888888889,
xtick={1,2,3,4,5,6},
xticklabels={{0.5},{0.6},{0.7},{0.8},{0.9},{1}},
xlabel style={font=\Large},
xlabel={$\text{p}_\text{s}$},
ymin=-10,
ymax=85,
ylabel style={font=\Large},
ylabel={Percentage difference in solver-time},
axis background/.style={fill=white},
xmajorgrids,
ymajorgrids,
legend style={legend cell align=left, align=left, draw=white!15!black, font = \Large}
]
\addplot[ybar, bar width=0.178, fill=green, draw=black, area legend] table[row sep=crcr] {%
	1	63.2856228708297\\
	2	64.4795589955735\\
	3	63.4097708004983\\
	4	63.7534414095079\\
	5	62.9334953723679\\
	6	61.4108261088747\\
};
\addplot[forget plot, color=white!15!black] table[row sep=crcr] {%
	0.511111111111111	0\\
	6.48888888888889	0\\
};
\addlegendentry{$\gain_t = \zeros$}

\addplot[ybar, bar width=0.178, fill=blue, draw=black, area legend] table[row sep=crcr] {%
	1	39.9569573596766\\
	2	41.7148453856577\\
	3	40.765231276372\\
	4	40.7317648995701\\
	5	40.273194189859\\
	6	40.6464778110333\\
};
\addplot[forget plot, color=white!15!black] table[row sep=crcr] {%
	0.511111111111111	0\\
	6.48888888888889	0\\
};
\addlegendentry{$\gain_t = \tilde{\gain}_t$ as in \eqref{e:gain2}}

\addplot[ybar, bar width=0.178, fill=black, draw=black, area legend] table[row sep=crcr] {%
	1	-0.521730023213216\\
	2	1.95202819397225\\
	3	0.252368857459896\\
	4	-0.526874018307626\\
	5	-1.55057041696699\\
	6	1.60609227380304\\
};
\addplot[forget plot, color=white!15!black] table[row sep=crcr] {%
	0.511111111111111	0\\
	6.48888888888889	0\\
};
\addlegendentry{without \eqref{e:stability constraint}}

\end{axis}

\begin{axis}[%
width=5.902in,
height=4.322in,
at={(0in,0in)},
scale only axis,
xmin=0,
xmax=1,
ymin=0,
ymax=1,
axis line style={draw=none},
ticks=none,
axis x line*=bottom,
axis y line*=left,
legend style={legend cell align=left, align=left, draw=white!15!black}
]
\end{axis}
\end{tikzpicture}%
	\end{adjustbox}	
	\caption{Percentage difference in solver-time with respect to the Theorem \ref{th:stbl} when $p_s$ varies and $p_c = 0.8$, $\authority = 5$ remain fixed.}
	\label{fig:runtime_snoise}
\end{figure}

 We repeated the above experiments for one sample path on intel i7-8750, 6 cores, 12 threads processor with 16 GB DDR4 RAM without invoking parallel pool in MATLAB. In one sample path, there are 120 time steps and 40 optimization instants. For each optimization instant, we compute the percentage difference with respect to the case ``$\gain_t$ as in \eqref{e:gain}'' and then take average over 40 optimization instants. The solver-time for the case ``$\gain_t$ as in \eqref{e:gain}'' is considered the base-value for the comparison. Our observations are given below:
\begin{enumerate}[label={\rm (O\arabic*)},leftmargin = *, nosep, start = 7]
	\item In Fig.\ \ref{fig:runtim_authority}, we plot the percentage difference in the solver-time with respect to  $\authority$ picked from the set $\{2,3,4,5,10\}$ while $p_s = p_c = 0.8$ remain fixed as in the observation \ref{o:umsb}. When $\gain_t = \zeros$, the solver-time is around $64\%$ less than the base-value with slight variations when we vary $\authority$. When $\gain_t = \tilde{\gain}_t$, the solver-time is around $44\%$ less than the base-value with slight variations when we vary $\authority$. The absence of \eqref{e:stability constraint} does not affect much. \label{o:runtime_authority}
\item In Fig. \ref{fig:runtime_cnoise}, we plot the percentage difference in the solver-time with respect to $p_c$ picked from the set $\{i/10 \mid i = 5, \ldots 10 \}$ while $p_s = 0.8$, $\authority = 5$ remain fixed as in the observation \ref{o:cmsb}. In this experiment, the observations are same as in \ref{o:runtime_authority} with a very slight difference. 
\item In Fig. \ref{fig:runtime_snoise}, we plot the percentage difference in the solver-time with respect to $p_s$ picked from the set $\{i/10 \mid i = 5, \ldots 10 \}$ while $p_c = 0.8$, $\authority = 5$ remain fixed as in the observation \ref{o:smsb}. In this experiment also, the observations are same as in \ref{o:runtime_authority} with a very slight difference.
\end{enumerate}	  

\par
\begin{figure}[h]
	\centering
	\begin{adjustbox}{width=\columnwidth}
%		\input{mae_correlated}
\begin{tikzpicture}

\begin{axis}[%
width=3.028in,
height=2.754in,
at={(2.509in,1.111in)},
scale only axis,
xmin=1,
xmax=6,
xtick={1,2,3,4,5,6},
xticklabels={{0.5},{0.6},{0.7},{0.8},{0.9},{1}},
xlabel={$p_{gx}$},
xmajorgrids,
ymin=6.5,
ymax=8,
ylabel={Empirical MAE},
ymajorgrids,
axis background/.style={fill=white},
legend style={at={(0.223,0.762)},anchor=south west,legend cell align=left,align=left,draw=white!15!black}
]
\addplot [color=red,dashed,line width=2.0pt]
table[row sep=crcr]{%
	1	6.50501292751256\\
	2	6.95798183182455\\
	3	7.27560580456292\\
	%4	7.53053795997106\\
	4	7.45174276775894\\
	5	7.72131693889378\\
	6	7.91375362243305\\
};
\addlegendentry{$p_{gx} = p_{gc}$};

\addplot [color=green,solid,line width=2.0pt]
table[row sep=crcr]{%
	1	7.39761380741959\\
	2	7.40745971340694\\
	3	7.48766745837241\\
	4	7.45174276775894\\
	5	7.46674009032707\\
	6	7.47049998160288\\
};
\addlegendentry{$p_{gx} = p_{gs}$};

\end{axis}
\end{tikzpicture}%
	\end{adjustbox}
	\caption{Empirical MAE with Gilbert-Elliott channel model}
	\label{Fig:correlated_mae}
\end{figure}
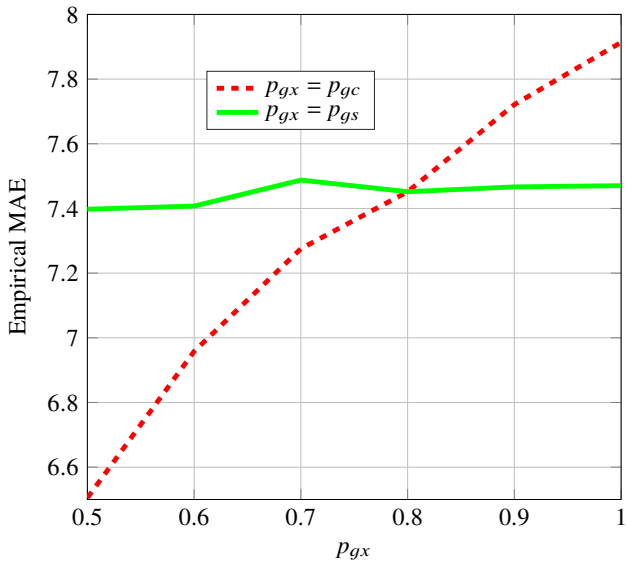

 In this article, we proposed an analytical framework, which is applicable to control with mutually independent and individually i.i.d. channels. For numerical experiments we can also consider that channels are mutually independent but packet dropouts in each channel are correlated. Such channels are often modelled by the Gilbert-Elliott channel model \cite{gilbert} as given in Fig.\ \ref{Fig:Gilbert-Elliott}. Here one assumes that each channel has two states -- good and bad, respectively. For simplicity, we model both channels similarly. For both channels, the bad state has successful transmission probability $p_b = 0$, the transition probabilities from good to bad state  and bad to good state are $p_{gb} = 0.2$ and $p_{bg} = 0.9$, respectively. The successful transmission probability of the good state of the control channel is $p_{gc}$ and that of the sensor channel is $p_{gs}$. Since both channels are modelled similarly except their good state, $p_{gx}$ in Fig.\ \ref{Fig:Gilbert-Elliott} takes value of $p_{gc}$ and $p_{gs}$, respectively, depending upon the channel. By considering the same experimental data as in the first part of this section, we did two experiments. In the first experiment, we fixed the good state of the sensor channel $p_{gx} = p_{gs} = 0.8$ and vary the good state of the control channel $p_{gx} = p_{gc}$ in the set $\{i/10 \mid i = 5, \ldots 10 \}$. In the second experiment, we fixed the good state of the control channel $p_{gx} = p_{gc} = 0.8$ and vary the good state of the control channel $p_{gx} = p_{gs}$ in the set $\{i/10 \mid i = 5, \ldots 10 \}$. We record the following observations:
	\begin{enumerate}[label={\rm (O\arabic*)},leftmargin = *, nosep, start = 10]
		\item Empirical MSB decreases with increase in $p_{gc}$ and fixed $p_{gs}$. It also decreases with increase in $p_{gs}$ and fixed $p_{gc}$. But in the second case the rate is lower than that in the first case; See Fig. \ref{Fig:correlated_msb}.  
		\item Empirical MAE per stage increases with increase in $p_{gc}$ and fixed $p_{gs}$. However, it does not vary much when $p_{gs}$ increases and $p_{gc}$ remain fixed; See Fig. \ref{Fig:correlated_mae}.  
	\end{enumerate}
The above numerical experiments show that our approach is computationally tractable in the presence of correlated channel noise as well. 

\begin{figure}[t]
	\centering
	\begin{adjustbox}{width=\columnwidth}
		\begin{tikzpicture}
		
		\begin{axis}[%
		width=3.028in,
		height=2.754in,
		at={(2.509in,1.111in)},
		scale only axis,
		xmin=1,
		xmax=6,
		xtick={1,2,3,4,5,6},
		xticklabels={{0.5},{0.6},{0.7},{0.8},{0.9},{1}},
		xlabel={$p_c$},
		xmajorgrids,
		ymin=6.5,
		ymax=9,
		ylabel={Empirical MAE},
		ymajorgrids,
		axis background/.style={fill=white},
		legend style={at={(0.023,0.694)},anchor=south west,legend cell align=left,align=left,draw=white!15!black}
		]
		\addplot [color=red,solid,line width=2pt]
		table[row sep=crcr]{%
			1	6.83383248645219\\
			2	7.26083749767367\\
			3	7.51274137456637\\
			4	7.75860955948248\\
			5	7.94149705520155\\
			6	8.11406794822533\\
		};
		\addlegendentry{$\gain_t$ as in \eqref{e:gain}};
		
		\addplot [color=green,dashed,line width=2pt]
		table[row sep=crcr]{%
			1	7.19024325017073\\
			2	7.71568836984241\\
			3	8.06397950925507\\
			4	8.3954607165089\\
			5	8.67364528204953\\
			6	8.94031219033195\\
		};
		\addlegendentry{$\gain_t = \zeros$};
		\addplot [color=blue,dotted,line width=2pt]
		table[row sep=crcr]{%
			1	6.99360104000791\\
			2	7.44242773827518\\
			3	7.72482464377026\\
			4	7.99986698930331\\
			5	8.23074295495548\\
			6	8.45049328595607\\
		};
		\addlegendentry{$\gain_t = \tilde{\gain}_t$ as in \eqref{e:gain2}};
		\addplot [color=black, mark = o]
		table[row sep=crcr]{%
			1	6.80583075034167\\
			2	7.22078073429472\\
			3	7.47078509703151\\
			4	7.70629743045316\\
			5	7.89355604438251\\
			6	8.0672607519323\\
		};
		\addlegendentry{without \eqref{e:stability constraint}}
						\end{axis}
		\end{tikzpicture}%
	\end{adjustbox}
	\caption{Empirical MAE per stage when $p_c$ varies from $0.5$ to $1$ while $\authority = 5$ and $p_s = 0.8$ remain fixed. }
	\label{Fig:varyingCnoiseEnergy}
\end{figure}
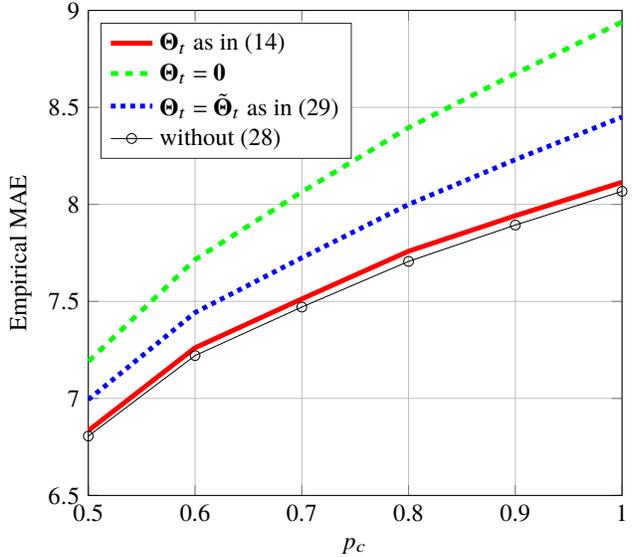

\begin{figure}[t]
	\centering
	\begin{adjustbox}{width=\columnwidth}
		\begin{tikzpicture}
		
		\begin{axis}[%
		width=3.028in,
		height=2.754in,
		at={(2.509in,1.111in)},
		scale only axis,
		xmin=1,
		xmax=6,
		xtick={1,2,3,4,5,6},
		xticklabels={{0.5},{0.6},{0.7},{0.8},{0.9},{1}},
		xlabel={$p_s$},
		xmajorgrids,
		ymin=7.5,
		ymax=8.5,
		ylabel={Empirical MAE},
		ymajorgrids,
		axis background/.style={fill=white},
		legend style={at={(0.423,0.54)},anchor=south west,legend cell align=left,align=left,draw=white!15!black}
		]
		\addplot [color=red,solid,line width=2pt]
		table[row sep=crcr]{%
			1	7.67140655340384\\
			2	7.70645333706253\\
			3	7.75856620582631\\
			4	7.75860955023601\\
			5	7.67544201871214\\
			6	7.632500980488\\
		};
		\addlegendentry{$\gain_t$ as in \eqref{e:gain}};
		
		\addplot [color=green,dashed,line width=2pt]
		table[row sep=crcr]{%
			1	8.00657948468713\\
			2	8.17603323628695\\
			3	8.31966159379826\\
			4	8.39546071663948\\
			5	8.38605314496341\\
			6	8.41503186931957\\
		};
		\addlegendentry{$\gain_t = \zeros$};
\addplot [color=blue,dotted,line width=2pt]
table[row sep=crcr]{%
	1	7.83763240901982\\
	2	7.9125268331308\\
	3	7.97845176671722\\
	4	7.99986698930331\\
	5	7.95578746950499\\
	6	7.95732951766673\\
};
\addlegendentry{$\gain_t = \tilde{\gain}_t$ as in \eqref{e:gain2}};	
\addplot [color=black, mark = o]
table[row sep=crcr]{%
	1	7.63313893998159\\
	2	7.66387028119047\\
	3	7.70950828116058\\
	4	7.70629743045316\\
	5	7.62942350283423\\
	6	7.58322917930202\\
};
\addlegendentry{without \eqref{e:stability constraint}}	
		\end{axis}
		\end{tikzpicture}%
	\end{adjustbox}
	\caption{Empirical MAE per stage when $p_s$ varies from $0.5$ to $1$ while $\authority = 5$ and $p_c = 0.8$ remain fixed. }
	\label{Fig:varyingSnoiseEnergy}
	
\end{figure}
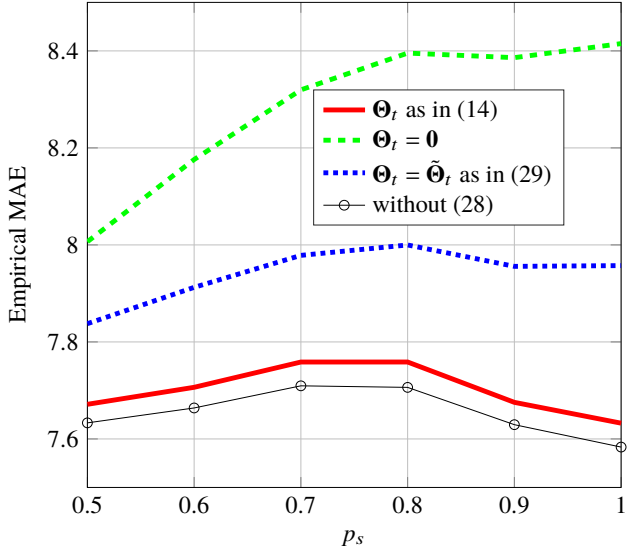

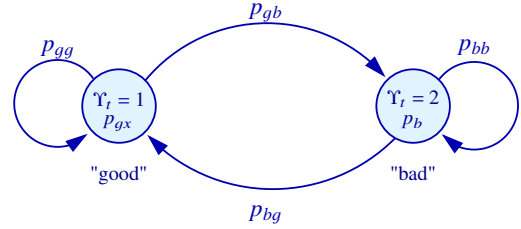
\begin{figure}[t]
	\centering
	\begin{adjustbox}{width=0.8\columnwidth}
		\begin{picture}(0,0)%
		\includegraphics{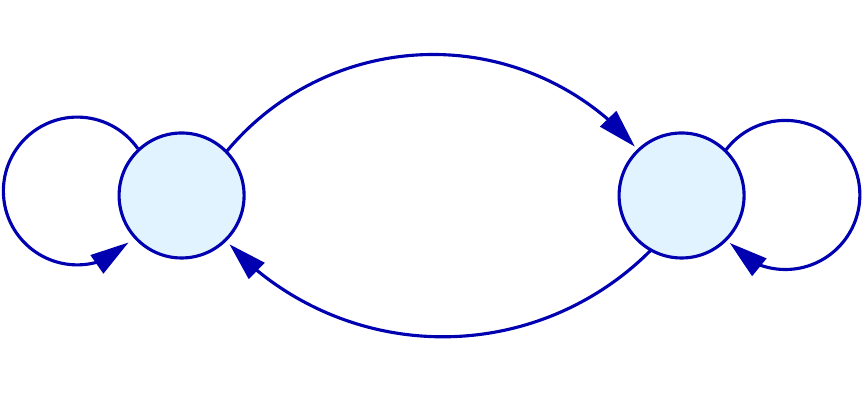}%
		\end{picture}%
		\setlength{\unitlength}{3947sp}%
		\begingroup\makeatletter\ifx\SetFigFont\undefined%
		\gdef\SetFigFont#1#2#3#4#5{%
			\reset@font\fontsize{#1}{#2pt}%
			\fontfamily{#3}\fontseries{#4}\fontshape{#5}%
			\selectfont}%
		\fi\endgroup%
		\begin{picture}(4142,1911)(-120,-1864)
		\put(3151,-1486){\makebox(0,0)[b]{\smash{{\SetFigFont{10}{12.0}{\familydefault}{\mddefault}{\updefault}{\color[rgb]{0,0,.56}"bad"}%
		}}}}
		\put(3651,-431){\makebox(0,0)[b]{\smash{{\SetFigFont{12}{14.4}{\familydefault}{\mddefault}{\updefault}{\color[rgb]{0,0,.69}$p_{bb}$}%
		}}}}
		\put(1951,-136){\makebox(0,0)[b]{\smash{{\SetFigFont{12}{14.4}{\familydefault}{\mddefault}{\updefault}{\color[rgb]{0,0,.69}$p_{gb}$}%
		}}}}
		\put(256,-436){\makebox(0,0)[b]{\smash{{\SetFigFont{12}{14.4}{\familydefault}{\mddefault}{\updefault}{\color[rgb]{0,0,.69}$p_{gg}$}%
		}}}}
		\put(1951,-1786){\makebox(0,0)[b]{\smash{{\SetFigFont{12}{14.4}{\familydefault}{\mddefault}{\updefault}{\color[rgb]{0,0,.69}$p_{bg}$}%
		}}}}
		\put(751,-1486){\makebox(0,0)[b]{\smash{{\SetFigFont{10}{12.0}{\familydefault}{\mddefault}{\updefault}{\color[rgb]{0,0,.56}"good"}%
		}}}}
		\put(747,-1036){\makebox(0,0)[b]{\smash{{\SetFigFont{10}{12.0}{\familydefault}{\mddefault}{\updefault}{\color[rgb]{0,0,.56}$p_{gx}$}%
		}}}}
		\put(751,-886){\makebox(0,0)[b]{\smash{{\SetFigFont{10}{12.0}{\familydefault}{\mddefault}{\updefault}{\color[rgb]{0,0,.56}$\Upsilon_t = 1$}%
		}}}}
		\put(3151,-867){\makebox(0,0)[b]{\smash{{\SetFigFont{10}{12.0}{\familydefault}{\mddefault}{\updefault}{\color[rgb]{0,0,.56}$\Upsilon_t = 2$}%
		}}}}
		\put(3151,-1036){\makebox(0,0)[b]{\smash{{\SetFigFont{10}{12.0}{\familydefault}{\mddefault}{\updefault}{\color[rgb]{0,0,.56}$p_b$}%
		}}}}
		\end{picture}%		
	\end{adjustbox}
	\caption{Transmission dropout model with a binary network state
		$(\Upsilon_t)_{t\in \{1,2\}}$: when $\Upsilon_t=1$ the channel is reliable with high successful transmission probabilities; $\Upsilon_t=2$ refers to a situation where the channel is
		unreliable and transmissions are more likely to be dropped.}
	\label{Fig:Gilbert-Elliott}
\end{figure}
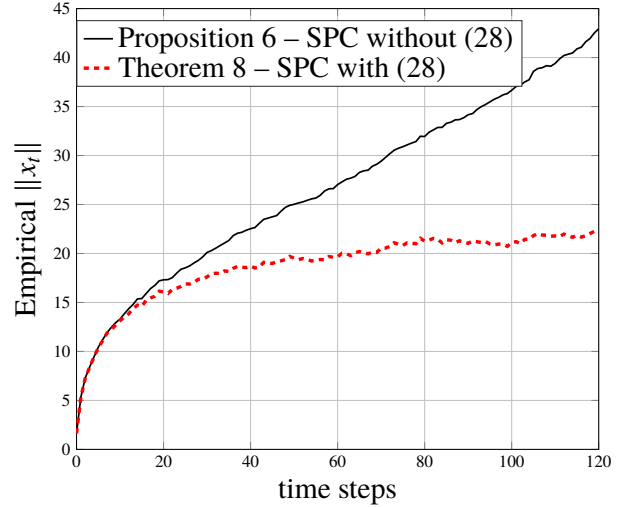
\begin{figure}
	\begin{adjustbox}{width = \columnwidth}
%		\input{yaverage_unstbl.tex}
\begin{tikzpicture}

\begin{axis}[%
width=4.221in,
height=3.566in,
at={(0.758in,0.481in)},
scale only axis,
xmin=0,
xmax=120,
xlabel={\Large{time steps}},
ymin=0,
ymax=45,
xmajorgrids,
ymajorgrids,
ylabel={\Large{Empirical $\norm{\st_t}$}},
axis background/.style={fill=white},
legend style={at={(0,0.823)},anchor=south west,legend cell align=left,align=left,draw=white!15!black}
]
\addplot [color=black,solid,line width=1.0pt]
table[row sep=crcr]{%
	0	1.64708388127294\\
	1	5.34014581733972\\
	2	7.29095490835059\\
	3	8.39536465384668\\
	4	9.42817999783708\\
	5	10.4278618485241\\
	6	11.1970639380426\\
	7	11.9115801327377\\
	8	12.4542290421198\\
	9	12.9396347437656\\
	10	13.2996607468379\\
	11	13.8556671915591\\
	12	14.3631887490681\\
	13	14.8000679213424\\
	14	15.3606424914283\\
	15	15.401394688388\\
	16	15.9135079005854\\
	17	16.4124626696417\\
	18	16.7439662629956\\
	19	17.2076070393881\\
	20	17.3108293702354\\
	21	17.3319946943137\\
	22	17.5282327169178\\
	23	17.9868640218931\\
	24	18.3918245153673\\
	25	18.5319384201006\\
	26	18.7365057441411\\
	27	19.0324342436326\\
	28	19.2966637892964\\
	29	19.6562094843156\\
	30	20.1078827425746\\
	31	20.2923804091012\\
	32	20.5421764720112\\
	33	20.810084571753\\
	34	21.066800418533\\
	35	21.3265707172764\\
	36	21.8223817557957\\
	37	22.0666598776407\\
	38	22.1630658485804\\
	39	22.3744180522911\\
	40	22.5304967742888\\
	41	22.6476653911025\\
	42	23.0745801829731\\
	43	23.4809720479622\\
	44	23.6093774576931\\
	45	23.7410712140789\\
	46	23.8717888950786\\
	47	24.3210857127292\\
	48	24.7085669684446\\
	49	24.9158250397088\\
	50	25.0226290806811\\
	51	25.1568579656024\\
	52	25.2659570668405\\
	53	25.4236460696206\\
	54	25.5713230717074\\
	55	25.6587829631197\\
	56	25.9373655801904\\
	57	26.3779192139201\\
	58	26.5997889355344\\
	59	26.6196066802742\\
	60	27.0346975161449\\
	61	27.3047495583281\\
	62	27.573834801938\\
	63	27.6873364688336\\
	64	27.8672372154587\\
	65	28.2699336116193\\
	66	28.4389602779486\\
	67	28.5017947191317\\
	68	28.9112295979321\\
	69	29.1117050928232\\
	70	29.4198419181529\\
	71	29.8015448761133\\
	72	30.1872110289622\\
	73	30.54187876869\\
	74	30.7459988225753\\
	75	30.9033334208088\\
	76	31.089688245258\\
	77	31.273194235199\\
	78	31.4245342507687\\
	79	31.965365156711\\
	80	31.9407087575656\\
	81	32.3727591399851\\
	82	32.6157401805019\\
	83	32.8661857102641\\
	84	32.8703328249765\\
	85	33.285300176574\\
	86	33.396760054588\\
	87	33.6339847370588\\
	88	33.6707847464764\\
	89	33.8638549468067\\
	90	34.165270102053\\
	91	34.2771781935732\\
	92	34.6845527192719\\
	93	34.9859011088615\\
	94	35.2116029534031\\
	95	35.4814841432474\\
	96	35.7461670460281\\
	97	35.9368695824608\\
	98	36.1930196960173\\
	99	36.3283154269636\\
	100	36.6617987561622\\
	101	37.0388047982844\\
	102	37.1920871806729\\
	103	37.509607828554\\
	104	37.7356163417409\\
	105	38.5852875435857\\
	106	38.882542318232\\
	107	38.9619868572949\\
	108	39.1684046104444\\
	109	39.1140491621994\\
	110	39.4114117398044\\
	111	39.8637692983555\\
	112	40.2122179185004\\
	113	40.3449764578333\\
	114	40.4680869992622\\
	115	40.8886584396371\\
	116	41.053746990091\\
	117	41.6078125341197\\
	118	41.9445285474526\\
	119	42.5146828889828\\
	120	42.9379353895654\\
};
\addlegendentry{\Large{Proposition \ref{th:main} -- SPC without \eqref{e:stability constraint}}};

\addplot [color=red,dashed,line width=2.0pt]
table[row sep=crcr]{%
	0	1.64708388127294\\
	1	5.3401511763019\\
	2	7.2533653729228\\
	3	8.36968169531703\\
	4	9.43242326257864\\
	5	10.3655693692909\\
	6	11.0915772624433\\
	7	11.8752831403702\\
	8	12.2332435522976\\
	9	12.6989998370346\\
	10	13.1481158320857\\
	11	13.5155975087955\\
	12	13.8491344740816\\
	13	14.3727925000451\\
	14	14.7263730630579\\
	15	14.5935869735497\\
	16	15.2276266028252\\
	17	15.4595724292253\\
	18	15.6411616536452\\
	19	16.1447643498898\\
	20	16.085705113976\\
	21	15.9241635214889\\
	22	16.3116164576548\\
	23	16.4483935579691\\
	24	16.6354895036874\\
	25	16.8921937169233\\
	26	16.883624741151\\
	27	17.0606677751052\\
	28	17.367791615735\\
	29	17.4135602928812\\
	30	17.6064076033158\\
	31	17.9474637993999\\
	32	17.9785349772916\\
	33	17.9649741848983\\
	34	18.2204055060244\\
	35	18.159969531708\\
	36	18.4828200367914\\
	37	18.684481751158\\
	38	18.5961503404945\\
	39	18.5608651391139\\
	40	18.7175117675415\\
	41	18.5173976007379\\
	42	18.6015229364514\\
	43	19.0943119513126\\
	44	19.0305159427767\\
	45	18.9848684860137\\
	46	19.1563780720443\\
	47	19.2477594059923\\
	48	19.3847268054433\\
	49	19.7000259483157\\
	50	19.5223002662163\\
	51	19.3783791028437\\
	52	19.485275132383\\
	53	19.3953239587323\\
	54	19.2357239345026\\
	55	19.3136688471839\\
	56	19.3851079045694\\
	57	19.361813418641\\
	58	19.690011991792\\
	59	19.6160108336458\\
	60	19.6632544529139\\
	61	20.018166244158\\
	62	19.9716558041635\\
	63	19.7833045357911\\
	64	20.0763322098597\\
	65	20.2011051392133\\
	66	20.0608480971303\\
	67	19.9708438738794\\
	68	20.1323513942134\\
	69	20.1129021344058\\
	70	20.5367881943349\\
	71	20.5569334956604\\
	72	20.7166689716805\\
	73	21.0975791939905\\
	74	20.9306995065731\\
	75	20.8256169702734\\
	76	21.0943395932436\\
	77	21.0431280350815\\
	78	20.9919998197065\\
	79	21.5514810350784\\
	80	21.2612817263921\\
	81	21.3603633248654\\
	82	21.5419858625747\\
	83	21.4531153598309\\
	84	21.077055736853\\
	85	21.4052444418137\\
	86	21.3564004395181\\
	87	21.100691441285\\
	88	21.2950224926853\\
	89	21.1574315500389\\
	90	21.0376104278126\\
	91	21.2752998071645\\
	92	21.2058535415885\\
	93	21.1187593506576\\
	94	21.2206124768598\\
	95	20.9697434474006\\
	96	20.8744705501107\\
	97	20.9601047394887\\
	98	20.9363406923575\\
	99	20.731713932492\\
	100	21.1158724027644\\
	101	21.2032708957172\\
	102	21.0187551659662\\
	103	21.4408114785179\\
	104	21.4660621446123\\
	105	21.8231351872821\\
	106	22.0247976296692\\
	107	21.8894727857847\\
	108	21.8598549604132\\
	109	21.7575373308504\\
	110	21.790042313361\\
	111	21.7236332554016\\
	112	21.9763948229408\\
	113	21.7887298999517\\
	114	21.6549700529476\\
	115	21.8452486385096\\
	116	21.6753936346\\
	117	21.8077145359839\\
	118	22.0328816262537\\
	119	22.2666979875725\\
	120	22.188621269583\\
};
\addlegendentry{\Large{Theorem \ref{th:stbl} -- SPC with \eqref{e:stability constraint}}};

\end{axis}

\begin{axis}[%
width=5.833in,
height=4.375in,
at={(0in,0in)},
scale only axis,
xmin=0,
xmax=1,
ymin=0,
ymax=1,
hide axis,
axis x line*=bottom,
axis y line*=left,
legend style={legend cell align=left,align=left,draw=white!15!black}
]
\end{axis}
\end{tikzpicture}%
	\end{adjustbox}	
\caption{Comparison between Proposition \ref{th:main} and Theorem \ref{th:stbl}}
\label{fig:unstbl}
\end{figure}

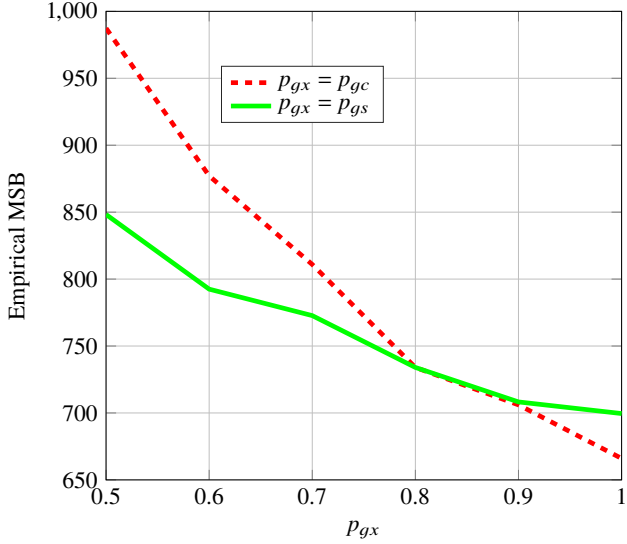
\begin{figure}[h]
	\centering
	\begin{adjustbox}{width=\columnwidth}
%		\input{msb_correlated}
\begin{tikzpicture}

\begin{axis}[%
width=3.028in,
height=2.754in,
at={(2.509in,1.111in)},
scale only axis,
xmin=1,
xmax=6,
xtick={1,2,3,4,5,6},
xticklabels={{0.5},{0.6},{0.7},{0.8},{0.9},{1}},
xlabel={$p_{gx}$},
xmajorgrids,
ymin=650,
ymax=1000,
ylabel={Empirical MSB},
ymajorgrids,
axis background/.style={fill=white},
legend style={at={(0.223,0.762)},anchor=south west,legend cell align=left,align=left,draw=white!15!black}
]
\addplot [color=red,dashed,line width=2.0pt]
table[row sep=crcr]{%
	1	987.406618854317\\
	2	877.093659426308\\
	3	810.917865860372\\
	%4	747.117602177358\\
	4   734.007471477405\\
	5	706.030688328673\\
	6	666.026899101558\\
};
\addlegendentry{$p_{gx} = p_{gc}$};

\addplot [color=green,solid,line width=2.0pt]
table[row sep=crcr]{%
	1	848.260897342835\\
	2	792.488856840493\\
	3	772.718637704645\\
	4	734.007471477405\\
	5	708.238335123799\\
	6	699.550081587088\\
};
\addlegendentry{$p_{gx} = p_{gs}$};

\end{axis}
\end{tikzpicture}%
	\end{adjustbox}
	\caption{Empirical MSB with Gilbert-Elliott channel model}
	\label{Fig:correlated_msb}
\end{figure}
We further consider a three dimensional stochastic LTI system \eqref{e:system} with matrices taken from \cite{ref:PDQ-15}:
\[ \A = \bmat{0&-0.80&-0.60\\ 0.80&-0.36&0.48\\ 0.60&0.48&-0.64}, \B = \bmat{0.16\\0.12\\0.14}, C=I_4, \] and simulation data $\stinit \sim N(0,I_4)$, $\wnoise_t \sim N(0,2I_3)$, $\mnoise_t \sim N(0,10I_4), p_s = p_c = 0.8, \authority = 15$. 
%\secondrevised{The above dynamical system can be written in the form of \eqref{e:decomposed st} with $\Aschur = 0.9$, $\Bschur = 0 $, \[\Aortho = \bmat{1 & 0 & 0\\0 & 0 & -1\\ 0 & 1 & 0} , \text{ and } \Bortho = \bmat{1 \\ 0 \\ 1}. \]
%	The reachability matrix (defined in \eqref{e:reachabilityMatrix})
%	\[ \reachab_{\reachindex}(\Aortho, \Bortho) = \bmat{
%		1 &	1&	1 \\
%		0&	-1&	0 \\
%		-1&	0	&1} 
%	\]
%	has full row rank for $\reachindex  = 3$.
%}
We repeatedly solved a finite-horizon CSOCP reported in Proposition \ref{th:main} and Theorem \ref{th:stbl} corresponding to states and control weights \[ Q = I_3 , Q_f = \bmat{12&-0.1&-0.4\\-0.1&19&-0.2\\-0.4&-0.2&2}, R =2 ,\] the optimization horizon \(N=4\), and compare their corresponding empirical $\norm{\st_t}$ in Fig. \ref{fig:unstbl}. This is improtant to note that the matrix $\A$ in this example is orthogonal and the matrix pair $(\A,\B)$ has reachability index $\reachindex  = 3$. Since the stability constraints \eqref{e:stability constraint} of Theorem \ref{th:stbl} suggest us to choose $N_r \geq \reachindex$, we choose $N_r = \reachindex = 3$ in our experiment. However, Proposition \ref{th:main} does not provide any such guideline to choose $N_r$, we take the standard choice $N_r = 1$ in our experiment. We simulated for $120$ time steps and took average of $500$ sample paths. We have the following observation:
\begin{enumerate}[label={\rm (O\arabic*)},leftmargin = *, nosep, start = 12]
	\item Empirical $\norm{\st_t}$ in case of Proposition \ref{th:main} increases almost linearly. However, in case of Theorem \ref{th:stbl} it is bounded below $23$ in $120$ time steps; See Fig. \ref{fig:unstbl}.    
\end{enumerate}
The above experiment demonstrates that there may exist a system, which can become unstable when stability constraints are not embedded in the corresponding optimization program under the given simulation data. 

%recalculation interval \(N_r = \reachindex = 3 \). 

\section{Epilogue}\label{s:epilogue}
We have employed a class of affine saturated received innovation feedback policies for a tractable formulation of the underlying CSOCP under the settings of unreliable channels. We proved that the proposed approach is tractable and ensures mean square boundedness of the controlled states. This approach can be extended for the class of policies parametrized in terms of dropouts along the lines of \cite{PDQ_Policy}. Moreover, inclusion of the joint chance constraints along the lines of \cite{paulson2015receding, RE-SPC, hokayem2013chance} and the sparsity in the control vector with help of a regularization term as in \cite{PDQ_Policy} will also be interesting extensions of the present work.  
% or


%\section{}
%
%\bibliographystyle{IEEEtran}  
\bibliographystyle{IEEEtran1}      
\bibliography{refs}          
\appendix
\section{Appendix}
	\begin{pf}[Proof of Lemma \ref{lem:ApproxEstimate}]
		On the event $\snoise_t = 1$, we can compute $\EE\left[ \st_t \mid \XX_t \right]$ as follows:
		\begin{align}\label{e:event1}
		& \EE \left[ \st_t \mid \XX_t \right] = \snoise_t \EE \left[ \st_t \mid \XX_t \right] \text{ since } \snoise_t = 1 \notag \\
		&= \snoise_t \EE \bigl[ \EE [\st_t \mid \YY_t, \XX_t ] \bigm|  \XX_t \bigr] = \snoise_t \EE \bigl[ \stfilt_t \bigm|  \XX_t \bigr] \notag \\
		&=  \EE \bigl[ \snoise_t \stfilt_t \bigm|  \XX_t \bigr] = \snoise_t \stfilt_t
		\end{align}
		On the event $\snoise_t = 0$, $\EE  \left[ \wnoise_{t-1} \mid \XX_t \right] = \zeros$ and $\EE  \left[\st_{t-1} \mid \XX_{t} \right] = \EE  \left[\st_{t-1} \mid \XX_{t-1} \right]$. Therefore, we can compute $\EE\left[ \st_t \mid \XX_t \right]$ as follows:
		\begin{align}\label{e:event0}
		& \EE\left[ \st_t \mid \XX_t \right] = (1-\snoise_t)\EE \left[\st_t \mid \XX_t \right] \text{ since } \snoise_t = 0 \notag \\
		&= (1-\snoise_t)\EE \left[\A\st_{t-1} + \B \control_{t-1}^a + \wnoise_{t-1} \mid \XX_t \right] \notag \\
		&= (1-\snoise_t)\EE  \left[\A\st_{t-1} + \B \control_{t-1}^a + \wnoise_{t-1} \mid \XX_{t} \right] \notag \\
		&= (1-\snoise_t)\left( \EE  \left[\A\st_{t-1} \mid \XX_{t} \right] + \B \control_{t-1}^a \right) \notag \\
		&= (1-\snoise_t)\left( \EE  \left[\A\st_{t-1} \mid \XX_{t-1} \right] + \B \control_{t-1}^a \right) \notag \\
		& = (1-\snoise_t)\left( \A\stest_{t-1} + \B \control_{t-1}^a \right) \text{ by definition of } \stest_t.
		\end{align}
		Since $\EE \left[ \st_t \mid \XX_t \right] = \snoise_t\left[ \st_t \mid \XX_t \right] + (1-\snoise_t)\EE \left[\st_t \mid \XX_t \right]$, the result follows by combining \eqref{e:event1} and \eqref{e:event0}. 
	\end{pf}
\begin{pf}[Proof of Lemma \ref{lem:objective}]
	Consider the objective function \eqref{e:cost_compact}.
	We substitute the stacked state vector \eqref{e:stacked state} in the objective function.
	% and stacked control vector \eqref{e:stacked control} . 
	\begin{align*}
	& V_t =  \EE_{\XX_t} \left[ \sum_{k = 0}^{N-1} (\norm{\st_{t+k}}^2_{Q_k} + \norm{\control_{t+k}}^2_{R_k}  ) + \norm{\st_{t+N}}^2_{Q_N}\right] \\
	& = \EE_{\XX_t} \left[ \norm{\calA\st_t + \calB\control_{t:N}^a + \calD \wnoise_{t:N}}^2_{\calQ} + \norm{\control_{t:N}^a}^2_{\calR} \right] \\
	& = \EE_{\XX_t} \Bigl[ \norm{\calA\st_t}^2_{\calQ} + \norm{\calD\wnoise_{t:N}}^2_{\calQ} + \norm{\control_{t:N}^a}^2_{\alpha} + 2(\st_t \transp \calA \transp \calQ \calB \\
	& \quad + \wnoise_{t:N} \transp \calD \transp \calQ \calB  ) \control_{t:N}  + 2\st_t\transp\calA\transp\calQ\calD\wnoise_{t:N} \Bigr].
	\end{align*}
	Let $\beta_t \Let \EE_{\XX_t}\left[ \norm{\calA\st_t}^2_{\calQ} + \norm{\calD\wnoise_{t:N}}^2_{\calQ} + 2\st_t\transp\calA\transp\calQ\calD\wnoise_{t:N} \right] = \EE_{\XX_t}\left[ \norm{\calA\st_t}^2_{\calQ} + \norm{\calD\wnoise_{t:N}}^2_{\calQ}  \right]$, then 
	\begin{equation} \label{e:objSolve}
	V_t = \EE_{\XX_t} \Bigl[\norm{\control_{t:N}^a}^2_{\alpha} + 2(\st_t \transp \calA \transp \calQ \calB + \wnoise_{t:N} \transp \calD \transp \calQ \calB  ) \control_{t:N}^a \Bigr] + \beta_t .
	\end{equation}
	We now substitute the stacked control vector \eqref{e:policyrepetitive out} in \eqref{e:objSolve} and for simplicity represent $\innovation_t = \meas_{t}-\hat{\meas}_{t} $ and $\noisyInnovation_{t} = \snoise_t \innovation_t $ to get the following set of equations:
	\begin{align*}
	& V_t = \EE_{\XX_t} \Bigl[\norm{\calG\offset_t + \calS \gain_t \ee(\noisyInnovation_{t:N})}^2_{\alpha} + 2(\st_t \transp \calA \transp \calQ \calB \\ 
	& \quad + \wnoise_{t:N} \transp \calD \transp \calQ \calB  ) \bigl( \calG\offset_t + \calS \gain_t \ee(\noisyInnovation_{t:N}) \bigr) \Bigr] + \beta_t\\
	& = \offset_t\transp \EE_{\XX_t} \Bigl[\calG\transp \alpha \calG \Bigr]\offset_t + \EE_{\XX_t} \Bigl[\norm{\calS\gain_t \ee(\noisyInnovation_{t:N})}^2_{\alpha}  \\
	& \quad + 2(\offset_t \transp \calG \transp \alpha + \st_t \transp \calA \transp \calQ \calB + \wnoise_{t:N} \transp \calD \transp \calQ \calB  ) \calS\gain_t \ee(\noisyInnovation_{t:N}) \Bigr] \\
	& \quad + \EE_{\XX_t} \Bigl[ 2(\st_t \transp \calA \transp \calQ \calB + \wnoise_{t:N} \transp \calD \transp \calQ \calB  )\calG\offset_t \Bigr] + \beta_t\\
	& = \offset_t\transp \Sigma_{\calG} \offset_t + \EE_{\XX_t} \Bigl[\norm{\calS\gain_t \ee(\noisyInnovation_{t:N})}^2_{\alpha} + 2\offset_t \transp \calG\transp\alpha\calS \gain_t \ee(\noisyInnovation_{t:N}) \\
	& \quad + 2(\st_t \transp \calA \transp \calQ \calB + \wnoise_{t:N} \transp \calD \transp \calQ \calB  ) \calS\gain_t \ee(\noisyInnovation_{t:N}) \Bigr] \\
	& \quad + 2\EE_{\XX_t} \Bigl[\st_t \transp \calA \transp \calQ \calB\calG\offset_t \Bigr] + \beta_t .
	\end{align*}
	Since $\EE_{\XX_t}\left[ \st_t \right] = \stest_t$, we obtain:
	\begin{equation}
	\begin{aligned}\label{e:Vt}
	& V_t = \offset_t\transp \Sigma_{\calG} \offset_t + \EE_{\XX_t} \Bigl[\norm{\calS\gain_t \ee(\noisyInnovation_{t:N})}^2_{\alpha} + 2(\offset_t \transp \calG\transp \alpha + \st_t \transp \calA \transp \calQ \calB \\
	& \quad + \wnoise_{t:N} \transp \calD \transp \calQ \calB  ) \calS \gain_t \ee(\noisyInnovation_{t:N}) \Bigr] + 2\stest_t \transp \calA \transp \calQ \calB\calG\offset_t + \beta_t .
	\end{aligned}
	\end{equation}
	Let us consider the term $\EE_{\XX_t} \left[ \offset_t \transp \calG \transp \alpha \calS \gain_t \ee(\noisyInnovation_{t:N})   \right]$ on the right hand side of \eqref{e:Vt} . By observing $\EE_{\XX_t} \left[ \ee_i( \noisyInnovation_{t+i}) \right] = \zeros$ for each $i = 1, \ldots, N$, we get
	\begin{align}\label{e:etay}
	\EE_{\XX_t} \left[ \offset_t \transp \calG\transp\alpha\calS \gain_t \ee(\noisyInnovation_{t:N})   \right] 
%	&= \snoise_t\offset_t \transp \Sigma_{\calG\calS}\bmat{\theta_{0,t} \\ \theta_{1,t}\\ \vdots \\ \theta_{N-1,t}}\ee_0(\meas_t - \hat{\meas}_t) \notag \\
	&=\offset_t \transp \Sigma_{\calG\calS} \gain_t^{(:,t)}\ee_0(\noisyInnovation_t),
	\end{align}
	where $\gain_t^{(:,t)} \Let \bmat{\theta_{0,t}\transp & \theta_{1,t} \transp & \ldots & \theta_{N-1,t} \transp}\transp $ represents the first $q$ columns of the gain matrix $\gain_t$. 
	Let us consider the term $\EE_{\XX_t} \Bigl[\norm{\calS \gain_t \ee(\noisyInnovation_{t:N})}^2_{\alpha} \Bigr]$ on the right hand side of \eqref{e:Vt}. In order to simplify offline computations, we perform the following manipulation:
	\begin{align}
	& \EE_{\XX_t} \Bigl[\norm{\calS\gain_t \ee(\noisyInnovation_{t:N})}^2_{\alpha} \Bigr] \nonumber \\
	& = \EE_{\XX_t} \Biggl[\norm{\calS\bmat{\gain_t^{(:,t)} & \gain_t^{\prime}} \bmat{\ee_0(\noisyInnovation_t) \\ \ee^{\prime}(\noisyInnovation_{t+1:N-1})} }^2_{\alpha} \Biggr] \nonumber \\
	& = \trace \Biggl( \Sigma_{\calS} \gain_t^{(:,t)} \ee_0(\noisyInnovation_t) \ee_0(\noisyInnovation_t) \transp (\gain_t^{(:,t)})\transp  \nonumber \\
	& \quad + \Sigma_{\calS} \gain_t^{\prime} \EE_{\XX_t}\Bigl[\ee^{\prime}(\noisyInnovation_{t+1:N-1})\ee^{\prime}(\noisyInnovation_{t+1:N-1})\transp\Bigr] (\gain_t^{\prime})\transp    \Biggr) \nonumber \\
	& = \trace(\Sigma_{\calS} \gain_t^{(:,t)} \Pi_{\meas_t} (\gain_t^{(:,t)})\transp) + \trace (\Sigma_{\calS} \gain_t^{\prime} \Sigma_{\ee} (\gain_t^{\prime})\transp) \label{e:Sigmay}, 
	\end{align}
	where $\Pi_{\meas_t} = \ee_0(\noisyInnovation_t) \ee_0(\noisyInnovation_t) \transp $ and $\Sigma_{\ee} = \\ \EE\Bigl[\ee^{\prime}(\noisyInnovation_{t+1:N-1})\ee^{\prime}(\noisyInnovation_{t+1:N-1})\transp \Bigr]$.
	We simplify the term $\EE_{\XX_t} \Bigl[ \wnoise_{t:N} \transp \calD \transp \calQ \calB \calS\gain_t \ee(\noisyInnovation_{t:N})\Bigr]$ in \eqref{e:Vt} as follows:
	\begin{align}
	& \EE_{\XX_t} \Bigl[\wnoise_{t:N} \transp \calD \transp \calQ \calB \calS \gain_t \ee(\noisyInnovation_{t:N})\Bigr] \nonumber \\
	&= \EE_{\XX_t} \Biggl[ \wnoise_{t:N} \transp \calD \transp \calQ \calB \calS \bmat{\gain_t^{(:,t)} & \gain_t^{\prime}} \bmat{\ee_0(\noisyInnovation_t) \\ \ee^{\prime}(\noisyInnovation_{t+1:N-1})}\Biggr] \nonumber \\
	&= \trace \left( \calD \transp \calQ \calB \mu_{\calS} \gain_t^{\prime} \EE_{\XX_t} \Biggl[\ee^{\prime}(\noisyInnovation_{t+1:N-1})\wnoise_{t:N} \transp \Biggr] \right) \nonumber \\
	&= \trace \left( \calD \transp \calQ \calB\mu_{\calS}\gain_t^{\prime} \Sigma_{\ee^{\prime}\wnoise} \right) \label{e:SigmaPsi},
	\end{align}
	where $\Sigma_{\ee^{\prime}\wnoise} = \EE \Biggl[\ee^{\prime}(\noisyInnovation_{t+1:N-1})\wnoise_{t:N} \transp \Biggr]$.
	Finally, we consider the term $\EE_{\XX_t} \Bigl[\st_t \transp \calA \transp \calQ \calB \calS\gain_t \ee(\noisyInnovation_{t:N})\Bigr]$ in \eqref{e:Vt} as follows:
	\begin{align}
	& \EE_{\XX_t} \Bigl[\st_t \transp \calA \transp \calQ \calB \calS \gain_t \ee(\noisyInnovation_{t:N})\Bigr] \notag \\
	& = \EE_{\XX_t} \Bigl[(\st_t - \stfilt_t)\transp \calA \transp \calQ \calB \calS \gain_t \ee(\noisyInnovation_{t:N})\Bigr] \notag \\
	& \quad + \EE_{\XX_t} \Bigl[ \stfilt_t\transp \calA \transp \calQ \calB \calS \gain_t \ee(\noisyInnovation_{t:N})\Bigr] \notag \\
	& = \EE_{\XX_t} \Bigl[e_t \transp \calA \transp \calQ \calB \calS \gain_t \ee(\noisyInnovation_{t:N})\Bigr] \notag \\ 
	& \quad +  \EE_{\XX_t} \Bigl[ \stfilt_t\transp \calA \transp \calQ \calB \calS \left(  \gain_t^{(:,t)} \ee_0(\noisyInnovation_{t}) + \gain_t^{\prime} \ee^{\prime}(\noisyInnovation_{t+1:N-1}) \right)\Bigr] \notag \\
	& = \EE_{\XX_t} \Bigl[e_t \transp \calA \transp \calQ \calB \calS \gain_t \ee(\noisyInnovation_{t:N})\Bigr] + \snoise_t \stest_t\transp \calA \transp \calQ \calB \mu_{\calS} \gain_t^{(:,t)} \ee_0(\innovation_t) \notag \\
	& =  \trace\Bigl(\calA \transp \calQ \calB \mu_{\calS} \gain_t^{\prime} \EE_{\XX_t} \Bigl[ \ee^{\prime} (\noisyInnovation_{t+1:N-1})e_t \transp \Bigr] \Bigr) \notag \\
	& \quad  + \stest_t\transp \calA \transp \calQ \calB \mu_{\calS} \gain_t^{(:,t)} \ee_0(\noisyInnovation_t) \notag \\
	& =  \trace\Bigl(\calA \transp \calQ \calB \mu_{\calS} \gain_t^{\prime} \Sigma_{e\ee^{\prime}} \Bigr)  + \stest_t\transp \calA \transp \calQ \calB \mu_{\calS} \gain_t^{(:,t)} \ee_0(\noisyInnovation_t), \label{e:Sigmae}
	\end{align}
	where $\Sigma_{e\ee^{\prime}} = \EE \Bigl[ \ee^{\prime} (\noisyInnovation_{t+1:N-1})e_t \transp \Bigr]$.
	Expression \eqref{e:obj out channel} follows by substituting \eqref{e:etay}, \eqref{e:Sigmay}, \eqref{e:SigmaPsi}, \eqref{e:Sigmae} in \eqref{e:Vt}, and ignoring the terms independent of the decision variables. 
	Therefore, the objective function in \eqref{e:cost_compact} is equivalent to \eqref{e:obj out channel} under the constraints \eqref{e:stacked state} and \eqref{e:policy out channel}.
\end{pf}
We present lemmas \ref{lem:msberror} - \ref{lem:stable_subsystem} before the proof of Lemma \ref{lem:ortho stable intermittent}.
	\begin{lem}\label{lem:msberror}
		Suppose that the estimator is driven by the recursion \eqref{e:est}, and let assumptions \ref{as:uncorrelated} -- \ref{as:lyapunov} hold. Then there exists $\tilde{\rho} > 0$ such that
		\begin{equation}\label{e:bound on dc}
		\EE\left[ \norm{\estError_{t}}^2 \mid \XX_{0}\right] \leq \tilde{\rho} \quad \text{ for all } t \geq 0 .
		\end{equation}
	\end{lem}
	\begin{pf}
		The estimation error is given by 
		\begin{equation}\label{e: estimation error}
		\estError_t = \stfilt_t - \stest_t = (1-\snoise_t)(\A \estError_{t-1} + \hat{\wnoise}_{t-1}).
		\end{equation}
		The matrix decomposition \eqref{e:decomposed st} allows us to write $\estError_t = \bmat{\estError_t^o \\ \estError_t^s}$ and $\hat{\wnoise}_t = \bmat{\hat{\wnoise}_t^o \\ \hat{\wnoise}_t^s}$, therefore
		\begin{align}\label{e:decomposed e_t}
		\estError_{t+1}^o &= (1-\snoise_{t+1})(\Aortho {\estError}_{t}^o + \hat{\wnoise}_t^o ) \\ 
		\estError_{t+1}^s &= (1-\snoise_{t+1})( \Aschur {\estError}_{t}^s + \hat{\wnoise}_t^s).
		\end{align}
		Let us first consider the expectation $\EE\left[ (\estError_{t+1}^s) \transp P \estError_{t+1}^s \mid \estError_{t}, \XX_0 \right]$ for some symmetric positive definite matrix $P \succ 0$. Since $\estError_t$ and $\hat{\wnoise}_t$ are independent at $t$, we get that
		\begin{align*}
		\EE\left[ (\estError_{t+1}^s) \transp P \estError_{t+1}^s \mid \estError_{t}, \XX_0 \right] & = (1-p_s) (\estError_{t}^s)\transp \Aschur \transp P \Aschur \estError_{t}^s \\
		& \quad + (1-p_s)\EE\left[(\hat{\wnoise}_t^s)\transp P \hat{\wnoise}_t^s \right]
		\end{align*}
		Since $\hat{\wnoise}_t$ has bounded variance for all $t \geq 0$. We can assume that there exists $C_1 < \infty$ such that $(1-p_s)\EE\left[(\hat{\wnoise}_t^s)\transp P \hat{\wnoise}_t^s \right] \leq C_1$ for all $t \geq 0$. Further, since $\Aschur$ is stable there exists some $\lambda_1 \in ]0,1[$ such that 
		\begin{align*}
		\EE\left[ (\estError_{t+1}^s) \transp P \estError_{t+1}^s \mid \estError_{t}, \XX_0 \right] \leq 
		(1-p_s)\lambda_1 \estError_{t}^s P \estError_{t}^s + C_1 \\
		\teL \lambda \estError_{t}^s P \estError_{t}^s + C_1
		\end{align*}
		where $\lambda = (1-p_s) \lambda_1 < 1$.  Since $\EE\left[ \norm{\estError_{t+1}^s}^2 \mid \estError_{t-1}, \XX_0 \right] = \EE \left[ \EE\left[ \norm{\estError_{t+1}^s}^2 \mid \estError_{t}, \XX_0 \right] \mid \estError_{t-1}, \XX_0 \right]$, we can show that
		\begin{align*}
		&	\EE\left[ (\estError_{t+1}^s) \transp P \estError_{t+1}^s \mid \estError_{0}, \XX_0 \right]  \leq \lambda^{t+1} (\estError_{0}^s) \transp P \estError_{0}^s + \dfrac{C_1}{1-\lambda}
		\end{align*}
		Therefore, 
		\begin{align*}
		\EE\left[ \norm{\estError_{t+1}^s}^2  \mid \XX_0 \right]  &\leq \lambda^{t+1}\frac{\lambda_{\max}(P)}{\lambda_{\min}(P)} \EE\left[ \norm{\estError_{0}^s}^2  \mid \XX_0 \right] \\
		& \quad + \dfrac{C_1}{(1-\lambda)\lambda_{\min}(P)},
		\end{align*}	
		where $\lambda_{\max}(P)$ and $\lambda_{\min}(P)$ are the largest and the smallest eigenvalues of $P$, respectively. 
		We can observe that on the event $\snoise_0 = 1$, $\estError_0 = 0$, whereas on the event $\snoise_0 = 0$, $\estError_0 = \stfilt_0$. Under the given initialization of the Kalman filter, we can easily see that $\stfilt_0 = K_0(\C \stinit + \mnoise_0) $. Since $\mnoise_0$ is independent of $\XX_0$ on the event $\snoise_0 = 0$, we can conclude that $\EE\left[ \norm{\estError_{0}^s}^2  \mid \XX_0 \right] = (1-\snoise_0)\EE\left[ \norm{\stfilt_0^s}^2 \right] \leq (1-\snoise_0)\trace(K_0\transp\Sigma_{\st_0}\C\transp)$. Therefore, there exists $C_2>0$ such that $\EE\left[ \norm{\estError_{t}^s}^2  \mid \XX_0 \right]  \leq C_2 $ for all $t \geq 0$. Similarly, we consider the expectation $\EE\left[ \norm{\estError_{t+1}^o}^2 \mid \estError_{t}, \XX_0 \right]$ and show that there exists $C_3 > 0$ such that $\EE\left[ \norm{\estError_{t}^o}^2  \mid \XX_0 \right]  \leq C_3 $ which implies $\EE\left[ \norm{\estError_{t}}^2  \mid \XX_0 \right] \leq C_2 + C_3 \teL \tilde{\rho}$ for all $t\geq 0$.    
	\end{pf}

\begin{lem} \label{lem: bound on finite hatwnoise}
	Consider \eqref{e:estEquation} and decomposition \eqref{e:decomposed st}. If there exist $m_3 > 0$ such that $\EE \left[ \norm{\hat{\wnoise}^o_t}^4 \right] \leq m_3$ for every $t$, then for some $h \in \N$, we have 
	\begin{equation}
	\EE \left[ \norm{\reachab_h(\Aortho, I_{\dortho})\hat{\wnoise}^o_{t  :h}}^4 \right] \leq m_3 h^4 \text{ for each } t.
	\end{equation}
\end{lem} 
\begin{pf}
	Let us begin with the term
	\begin{align*}
	& \norm{\reachab_h(\Aortho, I_{\dortho})\hat{\wnoise}^o_{t :h}}^4 = \left( \norm{\Aortho^{h-1} \hat{\wnoise}_t^o + \ldots + \hat{\wnoise}_{t+h-1}^o }  \right)^4  \\
	& \leq \left( \norm{ \hat{\wnoise}_t^o} + \ldots + \norm{\hat{\wnoise}_{t+h-1}^o }  \right)^4 \leq h^2 \left( \norm{ \hat{\wnoise}_t^o}^2 + \ldots + \norm{\hat{\wnoise}_{t+h-1}^o }^2  \right)^2 \\
	& \leq h^3 \left( \norm{ \hat{\wnoise}_t^o}^4 + \ldots + \norm{\hat{\wnoise}_{t+h-1}^o }^4  \right).
	\end{align*}
	We take expectation on both sides to get the bound
	\begin{align*}
	\EE \left[ \norm{\reachab_h(\Aortho, I_{\dortho})\hat{\wnoise}^o_{t :h}}^4 \right] & \leq h^3 \EE \left[ \norm{ \hat{\wnoise}_t^o}^4 + \ldots + \norm{\hat{\wnoise}_{t+h-1}^o }^4  \right] \\
	& \leq h^4 m_3 . 
	\end{align*}
\end{pf}
\begin{lem}\label{lem:fourth_moment_bound}
Let us consider $z_t$ as defined in Lemma \ref{lem:ortho stable intermittent}. There exists $m_2 > 0$ such that \[ \EE \left[ \norm{\tilde{\wnoise}_{\reachindex t + \ell}^o}^4  \mid z_0, \ldots , z_t \right] \leq m_2 \] for each $t$ and $\ell = 0, \ldots, \reachindex -1$.
\end{lem}
\begin{pf}
Let us consider the initialization $\estError_{-1} = \zeros$ and define 
	\begin{equation*}
	\tau_t \Let
	\begin{cases}
	\sup \{k \leq t \mid \snoise_k = 1 \} \quad & \text{ if }  \max_{k \leq t}{\snoise_k} = 1,\\
	-1 & \text{ otherwise },
	\end{cases}	 .
	\end{equation*}	
	then $\snoise_{\tau_t} = 1$ and from \eqref{e: estimation error} $\estError_{\tau_t} = 0$. In this way $ \estError_{\tau_t +1} =  \hat{\wnoise}_{\tau_t} 
	$
	and subsequently, $\estError_t = \reachab_{t-\tau_t}(\A, I) \hat{\wnoise}_{\tau_t: t - \tau_t}$ if $\tau_t \neq t$ and $\zeros$ otherwise. From \eqref{e:estEquation} we get
	\begin{equation}\label{e:estimated disturbance}
	\tilde{\wnoise}_t =
	\begin{cases}
	\reachab_{t-\tau_t + 1}(\A, I) \hat{\wnoise}_{\tau_t : t - \tau_t + 1} \quad & \text{ if }  \snoise_{t+1} = 1,\\
	\zeros & \text{ otherwise }.
	\end{cases}
	\end{equation}
	It is now clear that $\EE \left[ \norm{\tilde{\wnoise}_{\reachindex t + \ell}^o}^4  \mid z_0, \ldots , z_t , \snoise_0, \ldots , \snoise_{\reachindex t + \ell } \right] = \EE \left[ \norm{\tilde{\wnoise}_{\reachindex t + \ell}^o}^4  \mid \snoise_{\tau_{\reachindex t + \ell}}, \ldots , \snoise_{\reachindex t + \ell }\right]$. Let $\reachindex t + \ell + 1 - \tau_{\reachindex t + \ell} \Let h$,  then $ \tilde{\wnoise}_{\reachindex t + \ell}^o = \snoise_{\reachindex t + \ell}\reachab_h(\Aortho, I_{\dortho})\hat{\wnoise}^o_{\tau_{\reachindex t + \ell}:h}$. Since $\hat{\wnoise}_t$ is Gaussian for each $t$, there exists $m_3 > 0$ such that $\EE \left[ \norm{\hat{\wnoise}^o_{t}}^4  \right] \leq m_3$ for each $t$. We can conclude that        
\begin{align*}
& \EE \left[ \norm{\tilde{\wnoise}_{\reachindex t + \ell}^o}^4  \mid z_0, \ldots , z_t \right] \\ 
	&=p_s\sum_{h=1}^{\infty} \EE \left[ \norm{\reachab_h(\Aortho, I_{\dortho})\hat{\wnoise}^o_{\tau_{\reachindex t + \ell:h}}}^4  \mid \snoise_{\reachindex t + \ell - h+1},
	\ldots, \snoise_{\reachindex t + \ell} \right] \\
	& \quad p(\snoise_{\reachindex t + \ell - h+1} =1, \snoise_{\reachindex t + \ell - h+2} = \ldots = \snoise_{\reachindex t + \ell} =0) \\
	& \leq p_s \sum_{h=1}^{\infty}h^4 m_3 (1-p_s)^{h-1}p_s = m_3p_s^2\sum_{h=1}^{\infty}h^4  (1-p_s)^{h-1} \Let m_2. 
\end{align*}
	The last bound is implied by Lemma \ref{lem: bound on finite hatwnoise}.
\end{pf}
%=============================================================================
\begin{lem}\label{lem:sub-sampled_implies}
	Consider the system \eqref{e:decomposed st}. If there exists $\gamma_1 > 0$ such that $\EE_{\XX_0}\left[\norm{\tilde{\st}_{\reachindex t}^o}^2\right] \leq \gamma_1$ for all $t$, then there exists $\gamma_o > 0$ such that $\EE_{\XX_0}\left[\norm{\tilde{\st}_{t}^o}^2\right] \leq \gamma_o$ for all $t$.
\end{lem}
\begin{pf}
	This result is standard in literature, i.e. \cite{ref:RamChaMilHokLyg-10}. We provide proof for completeness. Letting $\norm{M}_F$ denote the Frobenius norm of a given matrix $M$, we compute a uniform bound on $\EE\left[\norm{\tilde{\st}_{\reachindex t + \ell}^o}^2\right]$ for $\ell = 0, \ldots, \reachindex -1$ as follows:
	\begin{align*}
	&\tilde{\st}_{\reachindex t + \ell}^o = \Aortho^{\ell}\tilde{\st}_{\reachindex t}^o + \reachab_{\ell}(\Aortho, \Bortho)\control_{\reachindex t: \ell}^a + \reachab_{\ell}(\Aortho, I)\tilde{\wnoise}_{\reachindex t:\ell}^o \\
	& \norm{\tilde{\st}_{\reachindex t + \ell}^o}^2 \leq 3 \left( \norm{\tilde{\st}_{\reachindex t}^o}^2 + \norm{\reachab_{\ell}(\Aortho, \Bortho)\control_{\reachindex t: \ell}^a}^2 + \norm{\reachab_{\ell}(\Aortho, I)\tilde{\wnoise}_{\reachindex t:\ell}^o}^2\right) 
	\end{align*} 
	%%%%%%%%%%%%%%%%%%%%
	% \trace(\Bortho \transp \Bortho) = \sum_{i=1}^m \lambda_i(\Bortho \transp \Bortho) \leq m \sigma_1(\Bortho)^2
	%%%%%%%%%%%%%%%%%%%%
	Since $\norm{\reachab_{\ell}(\Aortho, \Bortho)\control_{\reachindex t: \ell}^a}^2 \leq \norm{\reachab_{\ell}(\Aortho, \Bortho)}^2\norm{\control_{\reachindex t: \ell}^a}^2 \\ \leq \norm{\reachab_{\ell}(\Aortho, \Bortho)}_F^2\norm{\control_{\reachindex t: \ell}^a}_{\infty}^2\ell m \leq \norm{\reachab_{\ell}(\Aortho, \Bortho)}_F^2\authority^2\ell m\\
	= (\sum_{i=0}^{\ell-1}\norm{\Aortho^{i}\Bortho}_F^2)\authority^2\ell m = (\sum_{i=0}^{\ell-1}\trace(\Bortho\transp\Bortho))\authority^2\ell m \\
	\leq (\sum_{i=0}^{\ell-1}m \left(\sigma_1(\Bortho)\right)^2)\authority^2\ell m = \ell m \left(\sigma_1(\Bortho)\right)^2\authority^2\ell m$, we get 
	\[
	\norm{\tilde{\st}_{\reachindex t + \ell}^o}^2 \leq 3 \left( \norm{\tilde{\st}_{\reachindex t}^o}^2 + \left(\ell m \sigma_1(\Bortho)\authority \right)^2 + \ell \sum_{i=0}^{\ell-1}\norm{\tilde{\wnoise}_{\reachindex t +i}^o}^2 \right). 
	\]
	By taking the conditional expectation on both sides we get
	\begin{align*}
	& \EE_{\XX_0} \left[ \norm{\tilde{\st}_{\reachindex t + \ell}^o}^2 \right]  \leq 3 \left( \gamma_1 + \left(\ell m \sigma_1(\Bortho)\authority \right)^2\right) \\
	 & \quad + 3\ell^2 \max_{i=0, \ldots \ell-1} \EE_{\XX_0}\left[ \norm{\tilde{\wnoise}_{\reachindex t + i}^o}^2 \right] \\
	  & \quad \leq 3 \left( \gamma_1 + \left(\ell m \sigma_1(\Bortho)\authority \right)^2 + \ell^2 \sqrt{m_2}\right), 
	\end{align*}
	where $m_2$ is computed in Lemma \ref{lem:fourth_moment_bound}. 
	Therefore for $\ell = 0, \ldots,  \reachindex -1$ and for each $t$,
	\begin{align*}
	& \EE_{\XX_0}\left[ \norm{\tilde{\st}_{\reachindex t + \ell}^o}^2 \right] \leq 3 \left( \gamma_1 + (\reachindex m \sigma_1(\Bortho)\authority)^2 + \reachindex^2 \sqrt{m_2} \right)\teL \gamma_o.
	\end{align*} 
\end{pf}
\begin{lem}\label{lem:skip_free}
Consider the recursion \eqref{e:estEquation} and expression \eqref{e:z_difference}. The constraint \eqref{e:drift general 3} is satisfied when $\norm{\control^a_t}_{\infty} \leq \authority$ for all $t$.  	
\end{lem}
\begin{pf}
	  Since $\control_t^a$ is uniformly bounded there exists $m_1 > 0$ such that $\norm{\reachab_{\reachindex}(\Aortho, \Bortho)\control_{\reachindex t: \reachindex}^a}_{\infty} \leq m_1$ for all $t$. Let us consider $z_t$ defined in \eqref{e:z_difference}. 
	  Further, $\abs{(z_{t+1})^{(j)} - (z_t)^{(j)}} \leq \norm{z_{t+1} -z_t}_{\infty}$, the right side of which is bounded by $\norm{\reachab_{\reachindex}(\Aortho, \Bortho)\control_{\reachindex t: \reachindex}^a + \reachab_{\reachindex}(\Aortho, I_{\dortho})\tilde{\wnoise}_{\reachindex t: \reachindex}^o }_{\infty}$ by the definition of $z_t$. Therefore, we have $\abs{(z_{t+1})^{(j)} - (z_t)^{(j)}} \leq m_1 + \norm{\reachab_{\reachindex}(\Aortho, I_{\dortho})\tilde{\wnoise}_{\reachindex t: \reachindex}^o }_{\infty}$ for each $j = 1, \ldots, \dortho$. Let us consider the left hand side of \eqref{e:drift general 3}:
	\begin{align*}
	& \EE \left[ \abs{\left(z_{t+1}\right)^{(j)} - \left(z_t\right)^{(j)} }^4 \mid z_0^{(j)}, \ldots , z_t^{(j)} \right] \\
	& \quad \leq \EE \left[ \left(m_1 + \norm{\reachab_{\reachindex }(\Aortho, I_{\dortho})\tilde{\wnoise}_{\reachindex t: \reachindex}^o} \right)^4  \mid z_0^{(j)}, \ldots , z_t^{(j)} \right] \\
	& \quad \leq 8\left(m_1^4 + \EE \left[ \norm{\reachab_{\reachindex }(\Aortho, I_{\dortho})\tilde{\wnoise}_{\reachindex t: \reachindex}^o}^4  \mid z_0^{(j)}, \ldots , z_t^{(j)} \right]\right).
	\end{align*}
	Similar to the Lemma \ref{lem: bound on finite hatwnoise}, we get
	\begin{align*}
	& \EE \left[ \norm{\reachab_{\reachindex }(\Aortho, I_{\dortho})\tilde{\wnoise}_{\reachindex t: \reachindex}^o}^4  \mid z_0^{(j)}, \ldots , z_t^{(j)} \right] \\
	\quad & \leq \reachindex^4 \max_{i=0, \ldots, \reachindex-1}\EE \left[ \norm{\tilde{\wnoise}_{\reachindex t+i}^o}^4  \mid z_0^{(j)}, \ldots , z_t^{(j)} \right]\\
	\quad & \leq \reachindex^4 m_2,
	\end{align*}
	where the last bound is obtained by Lemma \ref{lem:fourth_moment_bound}. Therefore,
	\begin{align*}
	 \EE \left[ \abs{\left(z_{t+1}\right)^{(j)} - \left(z_t\right)^{(j)} }^4 \mid z_0^{(j)}, \ldots , z_t^{(j)} \right] \leq 8\left(m_1^4 + \reachindex^4 m_2\right) \teL M.
	\end{align*}
\end{pf}
\begin{lem}\label{lem:additive_noise_estimator}
	For a discrete time dynamical system \eqref{e:estEquation}, 
	\[ \EE[\tilde{\wnoise}_{t:\reachindex} \mid \XX_t] = \zeros \text{ for each } t \text{ and } \reachindex .\]
\end{lem}
\begin{pf}
We recall the expression \eqref{e:estimated disturbance}
\begin{equation*}
\tilde{\wnoise}_t = \snoise_{t+1}\reachab_{t-\tau_t + 1}(\A, I) \hat{\wnoise}_{\tau_t : t - \tau_t + 1}.
\end{equation*}
By taking the conditional expectation, we get
\begin{align*} 
&\EE[\tilde{\wnoise}_t \mid \XX_t] = p_s \EE[\reachab_{t-\tau_t + 1}(\A, I) \hat{\wnoise}_{\tau_t : t - \tau_t + 1} \mid \XX_t]\\
\quad  &= p_s \EE\left[\reachab_{t-\tau_t + 1}(\A, I) \EE[\hat{\wnoise}_{\tau_t : t - \tau_t + 1}\mid \YY_t] \mid \XX_t\right]= \zeros.
\end{align*}
 Similarly, $\EE[\tilde{\wnoise}_{t+\ell} \mid \XX_{t}] = \EE \left[ \EE[\tilde{\wnoise}_{t+\ell} \mid \XX_{t+\ell}] \mid \XX_{t+\ell}\right] = \zeros$ for $\ell = 0, \ldots, \reachindex - 1$. This completes the proof. 
\end{pf}
\begin{lem}\label{lem:stable_subsystem}
	There exists $\gamma_s > 0$ such that
	\[
	\sup_{t \in \Nz} \EE_{\XX_0}\left[\norm{\tilde{\st}_t^s}^2 \right] \leq \gamma_s .
	\]
\end{lem}	
\begin{pf}
A part of this proof is standard in literature, i.e. \cite{ref:RamChaMilHokLyg-10}. Let $P \succ 0$ be a symmetric positive definite matrix. 
\begin{align*}
&\EE[(\tilde{\st}_{t+1}^s) \transp P \tilde{\st}_{t+1}^s \mid \stest_t, \XX_0 ] = \norm{\tilde{\st}_t^s}_{\Aschur \transp P \Aschur}^2 + 2 (\tilde{\st}_t^s)\transp \Aschur \transp P \Bschur \control_t^a \\
\quad &+ \norm{\control_t^a}^2_{\Bschur \transp P \Bschur} + \EE[(\tilde{\wnoise}_t^s)\transp P \tilde{\wnoise}_t^s].
\end{align*}
Since $\Aschur$ is Schur stable, there exists some $\lambda \in ]0,1[ $ such that $\Aschur \transp P \Aschur \leq \lambda P$. We can also assume that there exists $m_4 \geq 0$ such that $\EE[(\hat{\wnoise}_t^s)\transp P \hat{\wnoise}_t^s] \leq m_4$ for each $t$ because $\hat{\wnoise}_t$ is Gaussian. We recall the expression \eqref{e:estimated disturbance} and compute $\EE[(\tilde{\wnoise}_t^s)\transp P \tilde{\wnoise}_t^s]$ as follows:
\begin{align*}
 &\EE[(\tilde{\wnoise}_t^s)\transp P \tilde{\wnoise}_t^s] = p_s \EE\left[\norm{\reachab_{t-\tau_t + 1}(\Aschur, I) \hat{\wnoise}^s_{\tau_t : t - \tau_t + 1}}_P^2 \right] \\
 &= p_s \sum_{h=1}^{\infty}\EE\left[\norm{\reachab_{h}(\Aschur, I) \hat{\wnoise}^s_{\tau_t : h}}_P^2 \mid t - \tau_t +1 = h\right]\\
 & \quad \quad \times  p(t - \tau_t +1 = h)
\\
&= p_s \sum_{h=1}^{\infty}\EE\left[\norm{\reachab_{h}(\Aschur, I) \hat{\wnoise}^s_{\tau_t : h}}_P^2 \mid t - \tau_t +1 = h\right]p_s(1-p_s)^{h-1}
\\
& \leq p_s^2 \sum_{h=1}^{\infty} h \sum_{i=1}^h \EE\left[\norm{\Aschur^{i-1}\hat{\wnoise}^s_{\tau_t + h-1}}_P^2 \mid t - \tau_t +1 = h\right](1-p_s)^{h-1} \\
& \leq p_s^2 \sum_{h=1}^{\infty} h \sum_{i=1}^h \lambda^{i-1} \EE\left[\norm{\hat{\wnoise}^s_{\tau_t + h-1}}_P^2 \mid t - \tau_t +1 = h\right](1-p_s)^{h-1}\\
& \leq p_s^2 \sum_{h=1}^{\infty} h m_4(1-p_s)^{h-1} \sum_{i=1}^h \lambda^{i-1} \\
& \leq \frac{m_4 p_s^2}{1-\lambda} \sum_{h=1}^{\infty} h (1-p_s)^{h-1} \teL m_5.
\end{align*}
Further, for $\varepsilon < \frac{1-\lambda}{\lambda}$ the Peter-Paul inequality provides us the bound $2 (\tilde{\st}_t^s)\transp \Aschur \transp P \Bschur \control_t^a \leq \varepsilon \norm{\tilde{\st}_t^s}_{\Aschur \transp P \Aschur}^2 + \frac{1}{\varepsilon} \norm{\control_t^a}^2_{\Bschur \transp P \Bschur}$. Therefore, 
\begin{align*}
&\EE[(\tilde{\st}_{t+1}^s) \transp P \tilde{\st}_{t+1}^s \mid \stest_t, \XX_0] \leq (1+\varepsilon)\norm{\tilde{\st}_t^s}_{\Aschur \transp P \Aschur}^2 \\
 & \quad \quad + (1+\frac{1}{\varepsilon})\norm{\control_t^a}^2_{\Bschur \transp P \Bschur} + m_5 \\
&  \leq 
(1+\varepsilon)\lambda\norm{\tilde{\st}_t^s}_{P}^2 + (1+\frac{1}{\varepsilon})m \lambda_{\max}(\Bschur \transp P \Bschur) \authority^2 + m_5.
\end{align*}
Defining $m_6 \Let (1+\frac{1}{\varepsilon})m \lambda_{\max}(\Bschur \transp P \Bschur) \authority^2 + m_5$ and iterating the above expression we get 
\begin{align*}
& \EE_{\XX_0}\left[(\tilde{\st}_{t+1}^s) \transp P \tilde{\st}_{t+1}^s \right] \leq ((1+\varepsilon)\lambda)^{t+1}\EE_{\XX_0} \left[\norm{\tilde{\st}_t^s}_{P}^2 \right] + \frac{m_6}{1-(1+\varepsilon)\lambda} \\
& \leq ((1+\varepsilon)\lambda)^{t+1}\lambda_{\max}(P)\snoise_0\EE_{\XX_0}\left[\norm{\stfilt_t^s}^2\right] + \frac{m_6}{1-(1+\varepsilon)\lambda} 
\\
& \leq ((1+\varepsilon)\lambda)^{t+1}\lambda_{\max}(P)\snoise_0\trace(K_0\transp\Sigma_{\st_0}\C\transp) + \frac{m_6}{1-(1+\varepsilon)\lambda}
\end{align*}
Since $(1+\varepsilon)\lambda < 1$, there exists $\gamma_s > 0$ such that $ \EE_{\XX_0}\left[\norm{\tilde{\st}_t^s}^2 \right] \leq \gamma_s$ for all $t$. 
\end{pf}
\begin{pf}[Proof of Lemma \ref{lem:ortho stable intermittent}]
 The equations in \eqref{e:drift general} are derived from Theorem \ref{t:PemRos-99} by substituting $z_t^{(j)}$ in place of $X_t$ in Theorem \ref{t:PemRos-99}. Therefore, there exists $\gamma_1 > 0$ such that $\EE_{\XX_0}\left[\norm{z_t}^2 \right] = \EE_{\XX_0}\left[\norm{\stest_{\reachindex t}^o}^2 \right] \leq \gamma_1$ for all $t$. Then by Lemma \ref{lem:sub-sampled_implies} there exists $\gamma_o > 0$ such that $\EE_{\XX_0}\left[\norm{\stest_{t}^o}^2 \right] \leq \gamma_o$ for all $t$.  
 Defining $\bar{\gamma} \Let \gamma_o + \gamma_s$ and by using Lemma \ref{lem:stable_subsystem} we get $\EE_{\XX_0}\left[\norm{\stest_{t}}^2 \right] \leq \bar{\gamma}$ for all $t$. 
 Now consider \eqref{e:drift general} and substitute \eqref{e:z_difference} to get that
 \begin{align*}
 & \EE_{\XX_{\reachindex t}} \left[ \left(z_{t+1} - z_t\right)^{(j)} \right]\\ 
 &= \EE_{\XX_{\reachindex t}} \left[ \left( (\Aortho^{\reachindex(t+1)})\transp\left( \reachab_{\reachindex}(\Aortho, \Bortho)\control_{\reachindex t: \reachindex}^a + \reachab_{\reachindex}(\Aortho, I_{\dortho})\tilde{\wnoise}_{\reachindex t: \reachindex}^o\right) \right)^{(j)}\right] \\
&= \EE_{\XX_{\reachindex t}} \left[ \left( (\Aortho^{\reachindex(t+1)})\transp \reachab_{\reachindex}(\Aortho, \Bortho)\control_{\reachindex t: \reachindex}^a  \right)^{(j)}\right], 
 \end{align*}
 where the last equality is due to Lemma \ref{lem:additive_noise_estimator}. Since \eqref{e:drift general 3} is satisfied by Lemma \ref{lem:skip_free}, \eqref{e:drift general} are equivalent to the following conditions for $t = 0, \reachindex, 2\reachindex, \ldots$, and $j = 1, \ldots, \dortho$:
 \begin{subequations} \label{e:drift intermittent intermediate}
 \begin{align}
 &	\Bigl( (\Aortho^{t+ \reachindex })\transp\reachab_{\reachindex}(\Aortho, \Bortho)\EE_{\XX_t} [\control_{t:\reachindex}^a] \Bigr)^{(j)} \leq -a \notag  \\
 & \quad \text{ whenever }  \left( (\Aortho^t) \transp \stest_{t}^o \right)^{(j)} > r, \label{e:drift1 intermittent intermediate} \\	
 &	\Bigl( (\Aortho^{t+  \reachindex})\transp\reachab_{\reachindex}(\Aortho, \Bortho)\EE_{\XX_t} [\control_{t:\reachindex}^a]\Bigr)^{(j)} \geq a \notag \\
 & 	\quad  \text{ whenever }  \left( (\Aortho^t) \transp \stest_{t}^o \right)^{(j)} < -r. \label{e:drift2 intermittent intermediate}
 \end{align}
 \end{subequations}
 Let us  define the component-wise saturation function  \(\R^{d_o}\ni z\longmapsto \sat_{r, \zeta}^\infty(z)\in\R^{d_o}\) to be 
 \[
 \bigl(\sat_{r, \zeta}^\infty(z)\bigr)^{(i)} = \begin{cases}
 z^{(i)} \zeta/r	& \text{if \(\abs{z^{(i)}} \le r\),}\\
 \zeta	 	& \text{if \(z^{(i)} > r\), and }\\
 -\zeta		& \text{otherwise,}
 \end{cases}
 \]
 for each \(i = 1, \ldots, d_o\). Now consider \eqref{e:kappa blocks control} with $\gain_t = \zeros$ and $(\offset_t)_{1:\reachindex m} = - \reachab_{\reachindex}(\Aortho, \Bortho)^\dagger \Aortho^{t+\reachindex} \sat_{r, \zeta}^\infty \bigl((\Aortho\transp)^{t} \stest_{ t}^o\bigr)$. It is clear that $\norm{\control_{t:\reachindex}}_{\infty} \leq \norm{\control_{t:\reachindex}}_2 \leq \sigma_1(\reachab_{\reachindex}(\Aortho, \Bortho)^\dagger)\sqrt{\dortho}\zeta \leq \authority$. Therefore, the given feedback policy satisfies \eqref{e:controlset}. The control sequence $\control_{ t:\reachindex}^a$ under the transmission protocol \ref{a:repetitive} is given by
 \begin{equation}\label{e:control sequence with tp3}
 \control_{t : \reachindex}^a = \tilde{\calG} \control_{ t:\reachindex},
 \end{equation} 
 where $\tilde{\calG} \in \R^{\reachindex m \times \reachindex m }$ is the principal submatrix of the diagonal matrix $\calG$ and $\EE[\tilde{\calG}^{(i,i)}] \geq p_c$ for each $i=1, \ldots, \reachindex m$. Therefore, the given policy also satisfies \eqref{e:drift intermittent intermediate} with $a=\zeta p_c$. We substitute $\EE[\control_{t:\reachindex}^a] \geq p_c \EE[\control_{t:\reachindex}]$ in \eqref{e:drift intermittent intermediate} to get \eqref{e:drift intermittent}. This completes the proof.  	
 \end{pf}

%=================================================================================
\begin{lem}
	\label{lem:estimator error bound cond}
Consider the system \eqref{e:system}. There exists $ \rho > 0 $ such that $ \EE_{\XX_0} \left[ \norm{e_t }^2 \right]  \leq \rho$ for all $t \geq 0$. 
\end{lem}
\begin{pf}
	Let us first observe that
	\[
	\EE_{\XX_0} \left[ \norm{e_t }^2 \right] =	\begin{cases}
	\EE_{\YY_0} \left[ \norm{e_t }^2 \right]	& \text{for } \snoise_0 = 1\\
	\EE_{\XX_0} \left[ \EE_{\YY_0} \left[ \norm{e_t }^2 \right] \right]		& \text{otherwise.}
	\end{cases}
	\]
	Since $\XX_0 \subset \XX_t$ for $t\geq 0$, the claim is implied by \cite[Lemma 4.2.2]{balakrishnan1987kalman} by using the tower property of the conditional expectation.
\end{pf}
\begin{pf} [Proof of Theorem \ref{th:stbl}]
	Since $\st_t = \st_t - \stfilt_t + \stfilt_t - \stest_t + \stest_t$, let us consider the inequality \[ \norm{\st_t}^2 \leq 3 \left( \norm{\st_t - \stfilt_t}^2 + \norm{\stfilt_t - \stest_t}^2 +  \norm{\stest_t}^2 \right). \]
	Let us apply the conditional expectation on the above inequality to obtain the bound: 
	\begin{align*}
	&\EE_{\XX_0}\left[ \norm{\st_t}^2 \right] \\
	& \quad \quad \leq 3\left( \EE_{\XX_0}\left[ \norm{\st_t - \stfilt_t}^2\right] + \EE_{\XX_0}\left[ \norm{\tilde{e}_t}^2\right] + \EE_{\XX_0}\left[ \norm{\stest_t}^2\right] \right)\\
	& \quad \quad \leq 3 \left( \rho + \EE_{\XX_0}\left[ \norm{\tilde{e}_t}^2\right] + \EE_{\XX_0}\left[ \norm{\stest_t}^2\right] \right) \text{ from Lemma \ref{lem:estimator error bound cond}} \\
	& \quad \quad \leq 3 \left( \rho + \tilde{\rho} + \EE_{\XX_0}\left[ \norm{\stest_t}^2\right] \right) \text{ from Lemma \ref{lem:msberror} }  \\
	& \quad \quad \leq 3 \left( \rho +  \tilde{\rho} +  \bar{\gamma} \right)  \text{ from Lemma \ref{lem:ortho stable intermittent}}\\
	& \quad \quad \teL \gamma \text{ for all } t \geq 0.
	\end{align*} 
\end{pf}
\end{document}